\documentclass[11pt,a4paper]{article}
\usepackage{amsmath, amsfonts, amssymb}
\usepackage{setspace, mathtools}
\usepackage{a4wide}
\usepackage{fullpage}
\usepackage{enumerate}

\usepackage{graphicx}
\usepackage{subfig}
\usepackage{caption}
\usepackage{sidecap}
\usepackage{wrapfig}

\usepackage{comment}

\usepackage[titletoc]{appendix}
\usepackage{color, subfig}
\usepackage{multirow, ulem}
\usepackage{amsthm}
\usepackage{caption}
\usepackage{tikz}
\usepackage{pgfplots}
\usetikzlibrary{calc}
\usepackage{tkz-euclide}
\usetikzlibrary{patterns}
\usepackage{hyperref}

\usepackage{tabu}

\hypersetup{
    colorlinks=true,
    linkcolor=blue,
    filecolor=magenta,      
    urlcolor=cyan,
    pdftitle={Sharelatex Example},
    bookmarks=true,
    citecolor=blue,
}

\usepackage{lipsum}

\usepackage[english]{babel}
\usepackage[nottoc]{tocbibind}
\usepackage{cite}
\usepackage{url}
\usepackage{authblk}

\usepackage{algorithm}
\usepackage{algpseudocode}

\theoremstyle{plain}
\newtheorem{definition}{Definition}

\newtheorem{theorem}{Theorem}
\newtheorem{proposition}{Proposition}
\newtheorem{corollary}{Corollary}

\title{\textbf{Facets of a mixed-integer bilinear covering set with bounds on variables}}
\author[1]{Hamidur Rahman}
\author[1]{Ashutosh Mahajan}
\affil[1]{\small Industrial Engineering and Operations Research, IIT Bombay, Mumbai, India. \protect\\ \texttt{{\{hrahman, amahajan\}@iitb.ac.in}}}

\date{}

\begin{document}

\maketitle

\begin{abstract}
We derive a closed form description of the convex hull of mixed-integer
bilinear covering set with bounds on the integer variables. This convex hull
description is determined by considering some orthogonal
disjunctive sets defined in a certain way. This description does not introduce
any new variables, but consists of exponentially many inequalities. An
extended formulation with a few extra variables and much smaller number of
constraints is presented. We also derive a linear time separation algorithm for
finding the facet defining inequalities of this convex hull. We study the
effectiveness of the new inequalities and the extended formulation using some
examples. 
\end{abstract}



\section{Introduction} \label{Sec1}

Consider the following mixed-integer bilinear covering set with bounds on the integer variables.
\[
S^U = \left\lbrace (x,y) \in \mathbb{Z}^n_+ \times \mathbb{R}^n_+ : \sum_{i = 1}^n x_i y_i \geq r, x \leq u \right\rbrace, \ \text{where } r > 0 \ \text{and } u \in \mathbb{N}^n \ \text{are given}.
\]
$S^U$ is a nonconvex set and even its continuous relaxation is nonconvex for $n \geq 2$. These constraints appear in the nonlinear formulation of the trim-loss problem (Harjunkoski et al. \cite{TrmLos}, Venderbeck \cite{vanderbeck}). In a trim-loss problem, we want to determine the best way to cut large rolls of raw materials into smaller pieces (or finals) using different patterns, so that the demand of finals is met. Let $N = \{ 1, \ldots, n \}$ be the index set that denotes the cutting patterns used, and $F$ be the index set of different sizes of the finals that are to be cut. Let $L$ be the size of each large roll and $l_j, j \in F$ be the lengths of the finals. The demands of the finals, say $d_j, j \in F$ are known. Let $x_{ij}$ be the number of final $j$ cut according in the pattern $i, i \in N, j \in F$, and $y_i$ be the number of rolls cut with cutting pattern $i, i \in N$. Therefore, we have the following constraints.

\begin{align}
& \sum_{i \in N} x_{ij} y_i \geq d_j, \ j \in F, \\
& \sum_{j \in F} l_j x_{ij} \leq L, i \in N. \label{TL}
\end{align}

Here, all the variables $x_{ij}$ and $y_i, i \in N, j \in F$ are non-negative integers. When the demands of the finals are high, we can consider the variables $y_i, i \in N$ as continuous variables without significantly affecting the optimal value. When integrality of $y$ can not be ignored then $S^U$ is a relaxation. Bounds on the variables $x_{ij}, i \in N, j \in F$ can be either given explicitly or be implicit from the knapsack constraints (\ref{TL}).

Let us consider a related set $S = \left\lbrace (x,y) \in \mathbb{Z}^n_+
\times \mathbb{R}^n_+ : \sum_{i = 1}^n x_i y_i \geq r \right\rbrace, r > 0$,
i.e., the set $S^U$ without the upper bounds on the variable $x$. Tawarmalani
et al. \cite{OrtDis} developed a scheme to get a tight convex relaxation using
orthogonal disjunctive subsets for a class of sets including $S$. They applied
the scheme to obtain the convex hull description of $S$ (denoted as
$conv(S)$) consisting of countably infinite number of facet defining
inequalities. But these facet defining inequalities of $conv(S)$ along with
the bound constraints are not sufficient to describe $conv \left( S^U
\right)$. Consider, for example,
\begin{align*}
\min & -x_1 + 10 y_1 - 2 x_2 + 12 y_2 \\
\text{s.t. } & x_1 y_1 + x_2 y_2 \geq 20, \\
& x_1 \leq 5, x_2 \leq 6, \tag{E} \label{E1} \\
& x_i \geq 0, y_i \geq 0, x_i \in \mathbb{Z}_+, i = 1, 2.
\end{align*}
Here, $r = 20, n = 2, u = (5,6)$. The point $(x_1, y_1, x_2, y_2) = \left( 5,
4, 6, 0 \right)$ is a global optimal solution with optimal value 23. But, if
we solve the relaxation defined by the facet defining inequalities of
$conv(S)$ (which we describe later), along with the bound constraints on $x$,
we get the solution $\omega = \left( 5, 1, 6, \frac{5}{6} \right)$ with
objective value 3. As expected, $\omega$ is not feasible for $S$ but lies in
$conv(S)$ because this point is the mid point of the two points $\left( 10, 2,
0, 0 \right), \left( 0, 0, 12, \frac{5}{3}\right) \in S$. Therefore, no facet
defining inequality of $conv(S)$ can cut off the point $\omega$ from $conv
\left( S^U \right)$. We will see later that the inequality $\frac{5 y_1}{20} +
\frac{6 y_2}{20} \geq 1$ is valid for $S^U$, and it cuts off the point $w$. In fact we show that this inequality is a facet defining inequality for $conv \left( S^U \right)$.

Optimizing a linear function over $S^U$ is a special case of nonconvex
(global) optimization problems. While nonconvex optimization problems are known to be generally
NP-Hard (Vavasis \cite{NP}, Katta and Santosh \cite{NPHard}), this particular
case turns out to be easily solvable. In almost all the algorithms for global
optimization, we take a convex relaxation of the feasible region and solve it
over successively refined partitions of the domain of the variables (Falk and
Soland \cite{FalkSol}, Belotti et al. \cite{IBM}, Sahinidis
\cite{Sahinidis2003}, Adjiman et al. \cite{AlphaBB}). A tighter relaxation
enables us to obtain tighter lower bounds on the problem and possibly
converges faster in a branch and bound framework.

There are different ways to get a convex relaxation depending on the function in a constraint. A bilinear function is a particular case of quadratic functions, for which there are several ways to get convex relaxations. McCormick relaxation (McCormick \cite{Mc}), Reformulation Linearization Technique (RLT) (Sherali \cite{RLT}, Sherali and Alameddin \cite{NewRLT}), Semidefinite relaxation (Anstreicher \cite{SDRLT}, \cite{ConAns}, Bao et al. \cite{REV}), Lagrangian relaxation (Voorhis \cite{Lag}) etc. are mostly used relaxation strategies of bilinear functions. Among these, McCormick and RLT give linear relaxations. However, these relaxations are generally weak in dimensions more than two.

The above mentioned relaxation strategies were developed for continuous
variables. These strategies can still be applied to get a convex relaxation
when some of the variables have integral restriction, but the relaxation is
generally even weaker. Furthermore, the above mentioned relaxation techniques
introduce new variables which naturally takes the problem to a higher
dimensional space. In order to obtain better bounds, one usually needs to exploit problem specific structures, like we do here.

In this article, we derive the closed form description of the convex hull of
the mixed-integer bilinear covering set $S^U$. We note that, the orthogonal
disjunctive technique of Tawarmalani et al. \cite{OrtDis} is not directly
applicable for the set $S^U$ to find  $conv \left( S^U \right)$. So, we relax
the orthogonal subsets of $S^U$ in such a way that the result is applicable.
Our work mainly addresses the following issues of the model of Tawarmalani et
al. Their model has infinitely many facet defining inequalities and these
inequalities along with the bound constraints gives us a weak relaxation of
our set. We show that $conv \left( S^U \right)$ is a polyhedron. We derive
both V-Polyhedron (i.e., description by sum of convex hull of the extreme
points and its recession cone) and H-Polyhedron (i.e., description by
intersection of finite number of half spaces) description of $conv \left( S^U
\right)$. We provide fast separation algorithms to find a violated facet defining inequality for both the sets $conv \left( S^U \right)$ and $conv(S)$. We also provide an extended formulation of $conv \left( S^U \right)$. Lastly, we provide some computational results that show the effectiveness of our cuts and the extended formulation.

Unless otherwise mentioned, we use the following notation throughout this article. For a given set $A$, we use $cl (A)$ to denote the closure of $A, conv(A)$ to denote the convex hull of $A, \mathcal{C}(A)$ to denote the conic hull of $A$  and $0^+(A)$ to denote the recession cone of $A$. $\mathbb{R}^n_+ = [0, \infty)^n = \{ x \in \mathbb{R}^n : x \geq 0 \}$. We use $N$ to denote the set $\{1, 2, \ldots, n \}$. For a point $(x,y) \in \mathbb{R}^n_+ \times \mathbb{R}^n_+$, we write $(x,y)$ in the form $(x_1, y_1, x_2, y_2, \ldots, x_n, y_n)$. We use $\mathcal{L}(i, x_i, y_i)$ to denote the point $(0,0, \ldots, x_i, y_i, \ldots, 0, 0)$, i.e., $x_j = 0, y_j = 0, \forall j \in N, j \neq i$.


\section{Convexification via Orthogonal Disjunction}

We start by a general result derived by Tawarmalani et al. \cite{OrtDis} for which some more notations are required. We use the same notation as in \cite{OrtDis} for convenience. Let $(z,u) \in \mathbb{R}^{\sum_{i=1}^n d_i} \times \mathbb{R}^{\sum_{i=1}^n d'_i}$, where $z_i \in \mathbb{R}^{d_i}$ and $u_i \in \mathbb{R}^{d'_i}$. Moreover, let us define the functions $t^j : \mathbb{R}^{\sum_{i=1}^n d_i} \times \mathbb{R}^{\sum_{i=1}^n d'_i} \rightarrow \mathbb{R}$ for $j \in J$, $v^k : \mathbb{R}^{\sum_{i=1}^n d_i} \times \mathbb{R}^{\sum_{i=1}^n d'_i} \rightarrow \mathbb{R}$ for $k \in K$ and $w^l : \mathbb{R}^{\sum_{i=1}^n d_i} \times \mathbb{R}^{\sum_{i=1}^n d'_i} \rightarrow \mathbb{R}$ for $l \in L$ where $J, K$ and $L$ are some index sets. Let us also define the sets $A \left( t^J,v^K,w^L \right)$ and $C \left( t^J,v^K,w^L \right)$ as below:
\begin{align*}
& A \left( t^J,v^K,w^L \right) = \left\lbrace (z,u) : t^j(z,u) \geq 1, \forall j \in J, v^k(z,u) \geq -1, \forall k \in K, w^l(z,u) \geq 0, \forall l \in L \right\rbrace, \text{and}\\
& C \left( t^J,v^K,w^L \right) = \left\lbrace (z,u) : t^j(z,u) \geq 0, \forall j \in J, v^k(z,u) \geq 0, \forall k \in K, w^l(z,u) \geq 0, \forall l \in L \right\rbrace.
\end{align*}

To describe the results, we need to additionally define positively-homogeneous functions. The following definition is taken from Rockafellar \cite{Rock} (1970).

\begin{definition}[Positively Homogeneous Function]
Let $f : \mathbb{R}^n \rightarrow [- \infty, \infty]$ be a function. $f$ is said to be a positively homogeneous function if, $f(\lambda x) = \lambda f(x), \forall \lambda > 0$.
\end{definition}

For example, $f(x,y) = \sqrt{x y}$ is positively homogeneous. Also, any linear function is positively homogeneous.

\begin{theorem}[Tawarmalani et al. \cite{OrtDis}] \label{ThMain}
Let $z=(z_1,...,z_i,...,z_n) \in \mathbb{R}^{\sum_{i=1}^n d_i}$, where $z_i \in \mathbb{R}^{d_i}$ and $Z \subseteq \mathbb{R}^{\sum_{i=1}^n d_i}$. Let $Z_i \subseteq Z$ for $i \in N = \{ 1,...,n \}$. Now let us consider the following assumptions:

\begin{description}

\item[A1:] $(z_1,...,z_i,...,z_n) \in Z_i \Rightarrow z_j = 0 \ \forall j \in N, j \neq i$,

\item[A2:] $conv(Z) = conv \left( \bigcup_{i=1}^n Z_i \right)$,

\item[A3:] $conv(Z_i) \subseteq proj_z(A_i) \subseteq cl(conv(Z_i))$, where
\[
A_i = \left\{ \mathcal{L}(i,z_i,u_i) : (z_i,u_i) \in A \left( t_i^{J_i},v_i^{K_i},w_i^{L_i} \right) \right\}
\]
such that $t_i^{j_i}, \ \forall j_i \in J_i, v_i^{k_i}, \ \forall k_i \in K_i$ and $w_i^{l_i}, \ \forall l_i \in L_i$ are positively-homogeneous functions for all $i \in N$ for some index sets $J_i, K_i$ and $L_i$, and $\mathcal{L}(i,z_i,u_i) = (0,...,0,z_i,u_i,0,...,0) \in \mathbb{R}^{\sum_{i=1}^n d_i} \times \mathbb{R}^{\sum_{i=1}^n d'_i}$,

\item[A4:] For all $i = 1,...,n$, $proj_z(C_i) \subseteq 0^+ \left( cl \left( conv \left( Z \right) \right) \right)$, where
\[
C_i = \left\{ \mathcal{L}(i,z_i,u_i) : (z_i,u_i) \in C \left( t_i^{J_i},v_i^{K_i},w_i^{L_i} \right) \right\},
\]
\end{description}
Then, $conv(Z) \subseteq proj_z(X) \subseteq cl(conv(Z))$, where,
\[
X = \left\{ (z,u) \left\vert   
            \begin{array}{l}
            
              \sum_{i=1}^n t_i^{j_i} (z_i,u_i) \geq 1, \hfill \forall (j_i)_{i \in N} \in \prod_{i=1}^n J_i, \\ \\
              
              \sum_{i \in I} v_i^{k_i} (z_i,u_i) \geq -1, \hfill \forall I \subseteq N, \forall (k_i)_{i \in I} \in \prod_{i \in I} K_i, \\ \\
              
              t_i^{j_i} (z_i,u_i) + v_i^{k_i} (z_i,u_i) \geq 0, \hfill \forall i \in N, \forall j_i \in J_i, \forall k_i \in K_i, \\ \\
              t_i^{j_i} (z_i,u_i) \geq 0, \forall i \in N, \hfill \forall j_i \in J_i, \\ \\
             
              w_i^{l_i} (z_i,u_i) \geq 0, \forall i \in N, \hfill \forall l_i \in L_i
              
             \end{array} \right.
             \right\}
\]
\end{theorem}

Using the above theorem, we can derive the convex hull for those sets which satisfy assumptions A1 - A4. Checking whether A1, A3 and A4 are satisfied by a given set is relatively easy. Verifying A2 might be difficult in practice. To overcome this difficulty, Tawarmalani et al. \cite{OrtDis} have used an alternative criterion, called convex extension property which is more general than assumption A2.



\begin{definition}[Convex Extension Property] \label{ConvExt}
Let $Z$ be a set in $\mathbb{R}^n$ and $Z_i \subseteq Z, i \in N$. The convex extension property holds for $Z$ if it satisfies the following two properties.
\begin{enumerate}[(i)]
\item If $z \in Z_i$, then $z_j = 0$ for all $j \in N, j \neq i$.
\item If $z \in Z$, then $z$ can be expressed as a sum of a convex combination of some points $\chi_i \in cl(conv(Z_i)), i \in N$ and a conic combination of rays $\psi_i \in 0^+(cl(conv(Z_i))), i \in N$, i.e. 
\[
z = \sum_{i \in N} \lambda_i \chi_i + \sum_{i \in N} \mu_i \psi_i   \tag{CE} \label{CE}
\]
\end{enumerate}
where $\mu_i \in \mathbb{R}_+, i \in N$ and $\lambda_i \in \mathbb{R}_+, i \in N$ with $\sum_{i \in N} \lambda_i = 1$.
\end{definition}

A collection of sets $Z_i, i \in N$ that satisfy condition (i) in Definition~\ref{ConvExt} are known as orthogonal sets. By definition a union of orthogonal sets satisfies the convex extension property. Some other sets that are not defined as union of orthogonal sets, for example, bilinear mixed-integer and pure-integer covering sets without variable bounds also satisfy this property. The convex extension property (\ref{CE}) is equivalent to the following criterion given in Tawarmalani et al. \cite{OrtDis}.

\[
cl (conv(Z)) = cl \left( conv \left( \bigcup_{i = 1}^n Z_i \right) \right) \tag{CE-P} \label{CE-P}
\]

Now, if we assume (\ref{CE}) or (\ref{CE-P}) instead of the assumption A2 in Theorem~\ref{ThMain}, we get $cl \left( proj_z X \right) = cl \left( conv \left( \bigcup_{i = 1}^n Z_i \right) \right) = cl(conv(Z))$ (Tawarmalani et al. \cite{OrtDis}). Since in many cases we only need $cl(conv(Z))$, it is useful to consider (\ref{CE}) or (\ref{CE-P}) instead of the assumption A2.


\section{On The Mixed-Integer Bilinear Covering Set $S$}

We start by revisiting the set $S = \left\lbrace (x,y) \in \mathbb{Z}^n_+ \times \mathbb{R}^n_+ : \sum_{i = 1}^n x_i y_i \geq r \right\rbrace, r > 0$, and the facet defining inequalities of its convex hull. Then we derive a property of extreme points of $conv(S)$ that we will later extend to $conv \left( S^U \right)$.

\subsection{The Convex Hull Description of $S$} \label{CvS}

Tawarmalani et al. \cite{OrtDis} showed that the set $S$ satisfies the assumptions A1, A3 and A4 of Theorem~\ref{ThMain} and the convex extension property (\ref{CE}) with respect to the orthogonal disjunctive subsets $S_i, i \in N$, where,
\[
S_i = \left\lbrace \mathcal{L}(i,x_i,y_i) \in \mathbb{Z}^n_+ \times \mathbb{R}^n_+ : x_i y_i \geq r \right\rbrace.
\]

Therefore, we can apply Theorem~\ref{ThMain} to construct the description of $conv(S)$. For this, first, we have to find the description of $conv \left( S_i \right)$. The continuous relaxation of the set $S_i$ is a convex set and the points of the form $\mathcal{L} \left( i,k,\frac{r}{k} \right), k \in \mathbb{N}$ are the extreme points of $conv \left( S_i \right)$. The convex hull description $conv(S_i)$ can be given as
\begin{equation}
conv \left( S_i \right) = \left\{ \mathcal{L}(i,x_i,y_i) : a_k x_i + b_k y_i \geq 1, k \in \mathbb{N} \right\}, \label{CSi}
\end{equation}
where $a_k x_i + b_k y_i = 1$ is the line passing through $\left(k, \frac{r}{k} \right)$ and $\left(k-1, \frac{r}{k-1} \right)$ for $k \in \mathbb{N} \setminus \{ 1 \}$ and $a_1 = 1, b_1 = 0$. Hence, we have $a_k = \frac{1}{2k - 1}$ and $b_k = \frac{k(k - 1)}{r(2k - 1)}$ for all $k \in \mathbb{N}$. We note that $conv(S_i)$ has countably infinite number of extreme points and facet defining inequalities. Consequently, $conv(S_i)$ is not a polyhedral set. Note that the recession cone $0^+(conv(S_i))$ of $conv \left( S_i \right)$ is the following set
\[
\left\{ (x,y) \in \mathbb{R}^n_+ \times \mathbb{R}^n_+ : x_j = 0, y_j = 0, j \in N, j \neq i \right\}.
\]

%
%
%
%

All sets $S_i, i \in N$ are identical to each other except for relabeling of indices. Thus, the coefficients $a_k$ and $b_k, k \in \mathbb{N}$ are identical for each $conv \left( S_i \right), i \in N$. Therefore, finding the coefficients $a_k, b_k, k \in \mathbb{N}$ for $conv \left( S_1 \right)$ is sufficient to get all the facets of $conv(S)$. The following collection of columns (\ref{M}) with countably infinite number of rows can be used to generate all the facet defining inequalities of $conv(S)$.

\[
\begin{bmatrix}
    x_1  & x_2 & x_3 & \dots & x_n \\
    a_2 x_1 + b_2 y_1  & a_2 x_2 + b_2 y_2 & a_2 x_3 + b_2 y_3 & \dots & a_2 x_n + b_2 y_n \\
    a_3 x_1 + b_3 y_1  & a_3 x_2 + b_3 y_2 & a_3 x_3 + b_3 y_3 & \dots & a_3 x_n + b_3 y_n \\ \tag{$M$} \label{M}
    \dots & \dots & \dots & \dots & \dots \\
    a_k x_1 + b_k y_1  & a_k x_2 + b_k y_2 & a_k x_3 + b_k y_3 & \dots & a_k x_n + b_k y_n \\
    \dots & \dots & \dots & \dots & \dots
\end{bmatrix}
\]

Theorem~\ref{ThMain} states that a facet defining inequality of $conv(S)$ is constructed by adding $n$ terms from (\ref{M}) taking exactly one term from each column and constraining their sum to be greater than or equal to one. All the facet defining inequalities are constructed this way. It is also clear that $conv(S)$ also has countably infinite number of facet defining inequalities. Since $a_k, b_k \geq 0, \forall k \in \mathbb{N}$, the recession cone $0^+(conv(S))$ of $conv(S)$ is the entire non-negative orthant $\mathbb{R}^n_+ \times \mathbb{R}^n_+$.

\subsection{Properties of The Extreme Points of $conv(S)$}

Here we derive the description of the extreme points of $conv(S)$ that we use later. We first note that $conv(S)$ is a closed set. This is because, if $(x,y) \notin conv(S)$, there exists a facet defining inequality of $conv(S)$ that strongly separates the point $(x,y)$ from $conv(S)$. Therefore, the point $(x,y)$ can not be a limit point of $conv(S)$, and consequently $conv(S)$ is a closed set. The convex extension property (\ref{CE-P}) applied to $S$ gives $conv(S) = conv \left( \bigcup_{i \in N} S_i \right)$.

\begin{theorem} \label{ExtThm}
$(\bar{x}, \bar{y})$ is an extreme point of $conv(S)$ if and only if $(\bar{x}, \bar{y})$ is an extreme point of $conv \left( S_i \right)$ for some $i \in N$.
\end{theorem}

\begin{proof}
Let $(\bar{x}, \bar{y})$ be an extreme point of $conv(S)$. If $(\bar{x}, \bar{y})$ belongs to $S_i$ for some $i \in N$, then it has to be an extreme point of $conv \left( S_i \right)$ as $S_i \subset S$. On the other hand, if $(\bar{x}, \bar{y})$ does not belong to any $S_i, i \in N$, then by convex extension property (\ref{CE-P}), $(\bar{x}, \bar{y})$ can be written as a convex combination of points in $S_i, i \in N$ which contradicts the extremality of the point $(\bar{x}, \bar{y})$.

Conversely, let $(\bar{x},\bar{y})$ be an extreme point of $conv \left( S_i \right)$ for some $i \in N$. Then, $\bar{x}_j = 0, \bar{y}_j = 0, \forall j \in N, j \neq i$. For contradiction, let $(\bar{x},\bar{y})$ be expressed as a convex combination of two distinct points $(\bar{x},\bar{y})^1$ and $(\bar{x},\bar{y})^2$ in $S$. Since $S \subset \mathbb{R}^n_+ \times \mathbb{R}^n_+$, then $\bar{x}^t_j = 0, \bar{y}^t_j = 0, \forall j \in N, j \neq i, t = 1, 2$. This implies that $(\bar{x},\bar{y})^1$ and $(\bar{x},\bar{y})^2$ belong to $S_i$. This is a contradiction to the fact that $(\bar{x},\bar{y})$ is an extreme point of $conv \left( S_i \right)$. Therefore, $(\bar{x},\bar{y})$ must be an extreme point of $conv(S)$.
\end{proof}

It is clear from Theorem~\ref{ExtThm} that any point of the form $\mathcal{L} \left(i, k, \frac{r}{k} \right), k \in \mathbb{N}$ is an extreme point of $conv(S)$ and vice versa, for all $i \in N$.


\section{On The Mixed-Integer Bilinear Covering Set $S^U$}

In this section we obtain a description of the convex hull of $S^U$ defined in
Section \ref{Sec1} and show that, unlike $conv(S)$, $conv \left( S^U \right)$ is
a polyhedron.

\begin{proposition} \label{Poly}
The set $conv \left( S^U \right)$ is a polyhedron.
\end{proposition}

\begin{proof}
Since there is an upper bound $u$ on the integer variable $x$, we have finitely many choices for $x$ in $S^U$. For each $i \in N$, we have $u_i + 1$ different choices for $x_i$. Since $x = 0$ is not a feasible choice for $S^U$, the total number of different choices for $x$ is $\prod_{i=1}^n (u_i + 1) - 1 = \eta$ (say). Let us denote them by $x^k, k =  1, \ldots, \eta$. Now, define the following polyhedral sets:
\[
F_k = \left\{ (x,y) \in \mathbb{Z}^n_+ \times \mathbb{R}^n_+ : \sum_{i=1}^n x_i y_i \geq r, x = x^k \right\}, k = 1, \ldots, \eta.
\]

Note that the set $F_k$ is constructed from $S^U$ by fixing $x = x^k$, and $S^U = \bigcup_{k = 1}^\eta F_k$. Further, the recession cone of $F_k$ is the set $\{ (0, y) \in \mathbb{R}^n \times \mathbb{R}^n : y \geq 0 \}$ for all $k = 1, \ldots, \eta$. Therefore, $S^U$ is a union of finite number of nonempty polyhedra with identical recession cones. So, from Corollary 4.44 in Conforti et al. \cite{CCZ}, we have $conv \left( S^U \right)$ is a polyhedron. 
\end{proof}

\subsection{The Extreme Point Description of $conv \left( S^U \right)$} \label{ExtDes}

%
%
%

Since $conv \left( S^U \right)$ is a polyhedron, it is closed and, therefore, it contains all its
extreme points. In this section we give a closed form description of the extreme points of $conv \left( S^U \right)$.

\begin{theorem} \label{ThExtDes}
Let $(\bar{x},\bar{y})$ be an extreme point of $conv \left( S^U \right)$. Then, $\bar{x}_t = p_t, \bar{y}_t = \frac{r}{p_t}$ for some $t \in N$, where $p_t \in \{ 1, \ldots, u_t \}$, and $\bar{x}_j \in \{ 0, u_j \}, \bar{y}_j = 0, \forall j \in N, j \neq t$, i.e., $(\bar{x},\bar{y})$ has the following form,
\[
\left( \bar{x}_1, 0, \bar{x}_2, 0, \ldots, \bar{x}_{t - 1}, 0, p_t,\frac{r}{p_t}, \bar{x}_{t + 1}, 0, \ldots, \bar{x}_n, 0 \right)
\]
where $p_t \in \{ 1, \ldots, u_t \}$ for some $t \in N, \bar{x}_j \in \{ 0, u_j \}, \forall j \in N, j \neq t$.
\end{theorem}

\begin{proof}
Let $(\bar{x}, \bar{y})$ be an extreme point of $conv \left( S^U \right)$.
Then $(\bar{x}, \bar{y}) \in S^U$. Therefore, $(\bar{x}, \bar{y}) \in F_k$ for
some $k \in \{1, \ldots, \eta \}$, and is an extreme point of $F_k$ where
$F_k$ is defined in the proof of Proposition~\ref{Poly}. Note that in the
description of $F^k$, there are $n$ bound constraints: $y_i \geq 0, i \in N$
and one linear constraint: $\sum_{i \in N} \bar{x}_i y_i \geq r$.

Since $(\bar{x}, \bar{y})$ is an extreme point of $F_k$, $n$ linear
constraints of $F_k$ must be active at $(\bar{x}, \bar{y})$. One can not
choose the $n$ constraints given by $y \geq 0$, otherwise $\sum_{i \in N}
\bar{x}_i y_i \geq r$ will be violated. So, the constraint $\sum_{i \in N}
\bar{x}_i y_i \geq r$ must be active. Therefore, there exists a $t \in N$ such that $\bar{x}_t \bar{y}_t = r$ and $y_j = 0, j \in N, j \neq t$.

We now show that $\bar{x}_j \in \{ 0, u_j \}$ for $\forall j \neq i$. If $\bar{x}_j \in (0, u_j), j \in N, j \neq i$, then $\bar{y}_j = 0$, and therefore, $(\bar{x}, \bar{y})$ can be written as a convex combination of the two points $(\bar{x}, \bar{y})^1$ and $(\bar{x}, \bar{y})^2$ having the exact same components as $(\bar{x}, \bar{y})$, except for the $j^{th}$ components of the variable $x$, and $\bar{x}^1_j = 0, \bar{x}^2_j = u_j$. Multipliers $1 - \lambda$ and $\lambda$ respectively provide the convex combination of $(\bar{x}, \bar{y})^1$ and $(\bar{x}, \bar{y})^2$, where $\lambda = \frac{\bar{x}_j}{u_j}$. This is a contradiction to the supposition that $(\bar{x}, \bar{y})$ is an extreme point of $conv \left( S^U \right)$.

Moreover, if $\bar{x}_j \in \{ 0, u_j \}$ for $j \neq i$, then we can not
write $(\bar{x}, \bar{y})$ as a convex combination of two different points in
$S^U$. This is because, if two such points exist, one of the points' $j^{th}$
component of the variable $x$ has to be more than $u_j$ or less than $0$,
neither of which is allowed. 
\end{proof}

\begin{corollary} \label{ExtNum}
$conv \left( S^U \right)$ has $2^{n - 1} \sum_{i = 1}^n u_i$ extreme points and $n$ extreme rays.
\end{corollary}

\begin{proof}
We see from the proof of Theorem~\ref{ThExtDes}, for a single choice of $\bar{x}_i \in \{ 1, \ldots, u_i \}$, we have $2^{n -1}$ different extreme points, and we have $\sum_{i = 1}^n u_i$ distinct such choices. Therefore, the total number of extreme points of $conv \left( S^U \right)$ is $2^{n - 1} \sum_{i = 1}^n u_i$, which is exponentially large, but finite. Consequently, $conv \left( S^U \right)$ is a polyhedral set.

On the other hand, we see that the recession cone $0^+ \left( conv \left( S^U \right) \right)$ of $conv \left( S^U \right)$ is the set $\{(x,y) \in \mathbb{R}^n_+ \times \mathbb{R}^n_+ : x = 0 \}$ which has $n$ extreme rays. 
\end{proof}

Note that Theorem~\ref{ThExtDes} and Corollary~\ref{ExtNum} give us the V-Description of $conv \left( S^U \right)$. We now turn our attention to the H-Description of $conv \left( S^U \right)$.

\subsection{The Convex Hull Description of $S^U$}

We have orthogonal disjunctive subsets of $S^U$,
\[
S^U_i = \left\{ \mathcal{L}(i,x_i,y_i) \in \mathbb{Z}^n_+ \times \mathbb{R}^n_+ : x_i y_i \geq r, x_i \leq u_i \right\}, i =1, \ldots, n.
\]
We note that $S^U_i \subset S^U$, and the recession cone of $cl \left( conv \left( S^U_i \right) \right)$ is the set:
\[
\left\{ (x,y) \in \mathbb{R}^n_+ \times \mathbb{R}^n_+ : x = 0, y_j = 0, \forall j \in N, j \neq i \right\}.
\]

We see that the assumption A1 of Theorem~\ref{ThMain} is satisfied by the set $S^U$ with respect to the orthogonal disjunctive subsets $S^U_i$. The polyhedral description of $conv \left( S^U_i \right)$ is
\[
conv \left( S^U_i \right) =  \left\{ \mathcal{L}(i, x_i, y_i) \in \mathbb{R}^n_+ \times \mathbb{R}^n_+ : a_k x_i + b_k y_i \geq 1, x_i \leq u_i, \forall k \in K_i  \right\},
\]
where $K_i = \{ 1, \ldots, u_i \}$, and as defined earlier, $a_k = \frac{1}{2k - 1}, b_k = \frac{k(k - 1)}{r(2k - 1)}, k \in K_i$. Therefore, assumption A3 of Theorem~\ref{ThMain} is satisfied by the set $S^U$ with respect to its orthogonal subsets $S^U_i, i \in N$.

On the other hand, the assumption A2 and convex extension property are not
satisfied by the set $S^U$ with respect to the subsets $S^U_i, i \in N$. An
extreme point of $conv \left( S^U \right)$ can have all $x$ components nonzero
which does not belong to $S^U_i$ for any $i \in N$. So, if it were in $conv \left( \bigcup_{i \in N} S^U_i \right)$, then it has to be a convex combination of two points in $S^U$ which contradicts the extremality of the point. 

In order to find the description of $conv \left( S^U \right)$, we use the following approach. The two inequalities $x_i y_i \geq r$ and $x_i \leq u_i$ in the description of $S^U_i$ together imply $y_i \geq \frac{r}{u_i}$. Let $\frac{r}{u_i} = \bar{u}_i$. Let us now define the following sets:
\[
S^L_i = \left\{ \mathcal{L}(i,x_i,y_i) \in \mathbb{Z}^n_+ \times \mathbb{R}^n_+ : x_i y_i \geq r, y_i \geq \bar{u}_i \right\}, i = 1, \ldots, n.
\]

By adding the lower bound on $y_i$ and ignoring the upper bound on $x_i$, we have a relaxation of $S^U_i$. The two sets $conv \left( S^U_i \right)$ and $conv \left( S^L_i \right)$ have exactly the same set of extreme points that are $u_i$ in number. Figure~\ref{FigComp1} and \ref{FigComp2} illustrate this observation.

\begin{figure}[h] 
\centering
\begin{minipage}{.5\textwidth}
  \centering
  \begin{tikzpicture}[scale=0.5]
    \draw[gray!50, thin, step=1] (0,0) grid (12,10);
    \draw[->] (0,0) -- (12.2,0) node[right] {$x_i$};
    \draw[->] (0,0) -- (0,10.2) node[above] {$y_i$};

   \foreach \x in {0,...,12} \draw (\x,0) -- (\x,0) node[below] {\scriptsize\x};
   \foreach \y in {0,...,10} \draw (0,\y) -- (0,\y) node[left] {\scriptsize\y};

   \fill[gray!50!,opacity=0.5] (1,10) -- (1,8) -- (2,4) -- (3,8/3) -- (4,2) -- (5,8/5) -- (6,8/6) --(6,10) -- cycle;
   \draw (1,10) -- (1,8) -- (2,4) -- (3,8/3) -- (4,2) -- (5,8/5) -- (6,8/6) --(6,10);
   
   \draw[fill] (1,8) circle [radius=2pt];
   \draw[fill] (2,4) circle [radius=2pt];
   \draw[fill] (3,8/3) circle [radius=2pt];
   \draw[fill] (4,2) circle [radius=2pt];
   \draw[fill] (5,8/5) circle [radius=2pt];
   \draw[fill] (6,8/6) circle [radius=2pt];
\end{tikzpicture}
\caption{$conv \left( S^U_i \right)$ for $r=8, x_i \leq 6$} \label{FigComp1}
\end{minipage}%
\begin{minipage}{.5\textwidth}
  \centering
  \begin{tikzpicture}[scale=0.5]
    \draw[gray!50, thin, step=1] (0,0) grid (12,10);
    \draw[->] (0,0) -- (12.2,0) node[right] {$x_i$};
    \draw[->] (0,0) -- (0,10.2) node[above] {$y_i$};

   \foreach \x in {0,...,12} \draw (\x,0) -- (\x,0) node[below] {\scriptsize\x};
   \foreach \y in {0,...,10} \draw (0,\y) -- (0,\y) node[left] {\scriptsize\y};

   \fill[gray!50!,opacity=0.5] (1,10) -- (1,8) -- (2,4) -- (3,8/3) -- (4,2) -- (5,8/5) -- (6,8/6) --(12,8/6) -- (12,10) -- cycle;
   \draw (1,10) -- (1,8) -- (2,4) -- (3,8/3) -- (4,2) -- (5,8/5) -- (6,8/6) --(12,8/6);
   
   \draw[fill] (1,8) circle [radius=2pt];
   \draw[fill] (2,4) circle [radius=2pt];
   \draw[fill] (3,8/3) circle [radius=2pt];
   \draw[fill] (4,2) circle [radius=2pt];
   \draw[fill] (5,8/5) circle [radius=2pt];
   \draw[fill] (6,8/6) circle [radius=2pt];
\end{tikzpicture}
\caption{$conv \left( S^L_i \right)$ for $r=8, y_i \geq \frac{8}{6}$} \label{FigComp2}
\end{minipage}
\end{figure}

We have the description of $conv \left( S^L_i \right)$ as following:
\begin{equation}
conv \left( S^L_i \right) = \left\{ \mathcal{L}(i, x_i, y_i) \in \mathbb{R}^n_+ \times \mathbb{R}^n_+ : a_k x_i + b_k y_i \geq 1, y_i \geq \bar{u}_i, \forall k \in K_i  \right\}, \label{CSLi0}
\end{equation}
where $K_i = \{ 1, \ldots, u_i \}$, and $a_k, b_k, k \in K_i$ are defined earlier. We also note that the recession cone of $conv \left( S^L_i \right)$ is the set:
\[
\left\{ (x,y) \in \mathbb{R}^n_+ \times \mathbb{R}^n_+ : x_j = 0, y_j = 0, j \in N, j \neq i \right\}.
\]
Let us now define a new set
\[
S^L = \bigcup_{i = 1}^n S^L_i.
\]

We will later derive the description of $conv \left( S^U \right)$ using $cl \left( conv \left( S^L \right) \right)$. We first observe that, since $S^L = \bigcup_{i=1}^n S^L_i$, we have,  
\[
cl \left( conv \left( S^L \right) \right) = cl \left( conv \left( \bigcup_{i=1}^n S^L_i \right) \right),
\]
i.e., the set $S^L$ satisfies the condition (\ref{CE-P}) with respect to the orthogonal disjunctive subsets $S^L_i, i \in N$.

\begin{proposition} \label{MyProp2}
The set $S^L$ satisfies all the assumptions A1 - A4 of Theorem~\ref{ThMain} with respect to the orthogonal disjunctive subsets $S^L_i, i \in N$.
\end{proposition}

\begin{proof}
We see that the assumption A1 holds from the definition of $S^L$. For the assumption A2, we have the convex extension property that is satisfied as noted above. Since we have the polyhedral description of $conv \left( S^L_i \right)$, the assumption A3 is satisfied. Lastly, we see that $0^+ \left( cl \left( conv \left( \bigcup_{i=1}^n S^L_i \right) \right) \right)$ is the entire non-negative orthant $\mathbb{R}^n_+ \times \mathbb{R}^n_+$, which implies that the assumption A4 is also satisfied. 
\end{proof}

We can now apply Theorem~\ref{ThMain} to obtain a description of $cl \left( conv \left( S^L \right) \right)$. We have $conv \left( S^L_i \right) = \left\{ \mathcal{L}(i, x_i, y_i) \in \mathbb{R}^n_+ \times \mathbb{R}^n_+ : a_k x_i + b_k y_i \geq 1, y_i \geq \bar{u}_i, k \in K_i \right\}$, where, $K_i = \{ 1, \ldots, u_i \}$ as defined earlier. Let us write it using a single index set as following:

\[
conv \left( S^L_i \right) = \left\lbrace \mathcal{L}(i,x_i,y_i) \in \mathbb{R}^n_+ \times \mathbb{R}^n_+ : l^{k_i} (x_i,y_i) \geq 1, k_i \in \bar{K}_i \right\rbrace,
\]
where, $\bar{K}_i = K_i \bigcup \{ u_i + 1 \}, l^{k_i}(x_i, y_i) = a_{k_i} x_i + b_{k_i} y_i$, where $a_{k_i} = \frac{1}{2{k_i} - 1}, b_{k_i} = \frac{k_i(k_i - 1)}{r(2k_i - 1)}, k_i \in K_i$ and $l^{(u_i + 1)_i}(x_i,y_i) = \frac{y_i}{\bar{u}_i}$. Note that the extreme points of $conv \left( S^L_i \right)$ are $\mathcal{L} \left( i, x_i, \frac{r}{x_i} \right), x_i = 1, \ldots, u_i$. Therefore, we have
\begin{equation}
conv \left( S^L_i \right) = conv \left( \left\{ \mathcal{L} \left( i, x_i,
\frac{r}{x_i} \right) : x_i = 1, \ldots, u_i \right\} \right) + \mathcal{C}
\left( \mathcal{L}(i, 1, 0), \mathcal{L}(i, 0, 1) \right), \label{CSLi}
\end{equation}
where $\mathcal{C} \left( \mathcal{L} \left( i, 1, 0 \right), \mathcal{L}(i, 0, 1) \right)$ is the conic hull of $ \left\{ \mathcal{L} \left( i, 1, 0 \right), \mathcal{L}(i, 0, 1) \right\}$. Now applying Theorem~\ref{ThMain} we have,
\begin{equation}
cl \left( conv \left( S^L \right) \right) = \left\lbrace (x,y) \in
\mathbb{R}^n_+ \times \mathbb{R}^n_+ : \sum_{i=1}^n l^{k_i} (x_i,y_i) \geq 1,
\forall \ (k_i)_{i=1}^n \in \prod_{i=1}^n \bar{K}_i \right\rbrace. \label{CSL}
\end{equation}

The set $cl \left( conv \left( S^L \right) \right)$ is a polyhedral set as it has finite number of facet defining inequalities in its description, and the number of facets is $\prod_{i = 1}^n |\bar{K}_i| = \prod_{i = 1}^n (u_i + 1)$ (which is exponentially large). Also, $0^+ \left( cl \left( conv \left( S^L \right) \right) \right)$ is the entire non-negative orthant $\mathbb{R}^n_+ \times \mathbb{R}^n_+$. Let us now derive some properties of the set $cl \left( conv \left( S^L \right) \right)$.

\begin{proposition}
The set $cl \left( conv \left( S^L \right) \right)$ is a polyhedral relaxation of $S^U$.
\end{proposition}

\begin{proof}
Since $S^L = \bigcup_{i = 1}^n S^L_i$, from (\ref{CSLi}) using Lemma~4.41 in \cite{CCZ} we have
\begin{equation}
cl \left( conv \left( S^L \right) \right) = conv \left( \bigcup_{i \in N} \left\{ \mathcal{L} \left( i, x_i \frac{r}{x_i} \right) : x_i = 1, \ldots, u_i \right\} \right) + \mathbb{R}^n_+ \times \mathbb{R}^n_+ \label{VCSL}
\end{equation}

Since $0^+ \left( conv \left( S^U \right) \right)$ is a subset of $\mathbb{R}^n_+ \times \mathbb{R}^n_+ = 0^+ \left( cl \left( conv \left( S^L \right) \right) \right)$, it is sufficient to show that all the extreme points of $conv \left( S^U \right)$ belong to $cl \left( conv \left( S^L \right) \right)$. Let $(\bar{x}, \bar{y})$ be an extreme point of $conv \left( S^U \right)$, then from Theorem~\ref{ThExtDes} we have
\[
(\bar{x}, \bar{y}) = \mathcal{L} \left( i, \bar{x}_i, \frac{r}{\bar{x}_i} \right) + (\bar{x}_1, 0, \ldots, \bar{x}_{i-1}, 0, 0, 0, \bar{x}_{i+1}, 0, \ldots, \bar{x}_n, 0)
\]
for some $i \in N$. This clearly shows that $(\bar{x}, \bar{y}) \in cl \left( conv \left( S^L \right) \right)$. 
\end{proof}

\begin{theorem} \label{MyTh1}
$(\bar{x},\bar{y})$ is an extreme point of $cl \left( conv \left( S^L \right) \right)$ if and only if $(\bar{x},\bar{y})$ is an extreme point of $conv \left( S^L_i \right)$ for some $i \in N$.
\end{theorem}

\begin{proof}
The statement follows from (\ref{VCSL}). 
\end{proof}

\begin{corollary}
$(\bar{x},\bar{y})$ is an extreme point of $cl \left( conv \left( S^L \right) \right)$ if and only if $(\bar{x},\bar{y})$ is an extreme point of $conv \left( S^U_i \right)$ for some $i \in N$.
\end{corollary}

\begin{proof}
Since $conv \left( S^U_i \right)$ and $conv \left( S^L_i \right)$ have exactly same set of extreme points, the result follows from Theorem~\ref{MyTh1}. 
\end{proof}

Here we observe that $cl \left( conv \left( S^L \right) \right)$ is a polyhedral relaxation of $S^U$ such that each extreme point of $cl \left( conv \left( S^L \right) \right)$ lies in $S^U$. Now we prove our main result.




\begin{theorem} \label{ThCH}
Let $\bar{S} = \{ (x, y) \in cl \left( conv \left( S^L \right) \right) : x \leq u \}$. Then, $conv \left( S^U \right) = \bar{S}$.
\end{theorem}

\begin{proof}
Minkowski Resolution Theorem (Theorem 4.15 in \cite{Bertsimas}) states that any polyhedral set having at least one
extreme point can be described by its extreme points and recession cone.  The polyhedral sets $conv \left( S^U \right)$ and $\bar{S}$
have the same recession cone $\{(x,y) \in \mathbb{R}^n_+ \times \mathbb{R}^n_+
: x = 0 \}$. 

The constraint $x_i \leq u_i$ passes through only one extreme point $\mathcal{L} \left( i, u_i, \frac{r}{u_i} \right)$ of
$cl \left( conv \left( S^L \right) \right)$ and does not cut off any of its extreme points. Therefore, adding
this constraint to $cl \left( conv \left( S^L \right) \right)$ only creates new extreme points of the form
\[
\left( w_1, 0, w_2, 0, \ldots, w_{i-1}, 0, p_i,\frac{r}{p_i}, w_{i+1}, 0,
\ldots, w_n, 0 \right),
\]
where $w_j \in \{ 0, u_j \}, j \in N, j \neq i, p_i \in \{ 1, \ldots, u_i \},
i \in N$. From Theorem~\ref{ThExtDes}, we see that such points lie in $S^U$, in fact, they are extreme points of $conv \left( S^U \right)$. Again, since $conv \left( S^U \right) \subseteq \bar{S}$, we have $\bar{S} = conv \left( S^U \right)$.
\end{proof}

\subsection{Facet Defining Inequalities of $conv \left( S^U \right)$}

We now focus our attention on the new inequalities that are generated by our procedure and their effectiveness. We have seen from Theorem~\ref{ThCH} that each facet defining inequality of $conv \left( S^U \right)$ is either a bound constraint $x_i \leq u_i$ for some $i \in N$ or a facet defining inequality of $cl \left( conv \left( S^L \right) \right)$ of the following form:
\[
\sum_{i=1}^n l^{k_i} (x_i,y_i) \geq 1, (k_i)_{i=1}^n \in \prod_{i=1}^n \bar{K}_i \tag{F\textsubscript{SL}} \label{FSL}
\]
where, $\bar{K}_i = K_i \bigcup \{ u_i + 1 \}, K_i = \{1, \ldots, u_i \}, l^{k_i}(x_i, y_i) = a_{k_i} x_i + b_{k_i} y_i, a_{k_i} = \frac{1}{2{k_i} - 1}, b_{k_i} = \frac{k_i(k_i - 1)}{r(2k_i - 1)}, k_i \in K_i$ and $l^{(u_i + 1)_i}(x_i,y_i) = \frac{y_i}{\bar{u}_i}, \bar{u}_i = \frac{r}{u_i}$. The inequality $\sum_{i=1}^n l^{k_i} (x_i,y_i) \geq 1$ is identical to one of the facet defining inequalities of $conv(S)$ if $(k_i)_{i=1}^n \in \prod_{i=1}^n K_i$. Now let $Q \subseteq N$ be a non-empty index set such that $k_i = (u_i + 1)_i$ for all $i \in Q$. Then, the inequalities of the form
\[
\sum_{i \in Q} \frac{y_i}{\bar{u}_i} + \sum_{i \in N \setminus Q} l^{k_i} (x_i,y_i) \geq 1, (k_i)_{i=1}^n \in \prod_{i=1}^n \bar{K}_i \tag{NF} \label{NF}
\]
are generated by applying our approach and they are not valid for $conv(S)$.


\subsection{An Extended Formulation of $conv \left( S^U \right)$} \label{SecExtFor}

We saw that the description of $conv \left( S^U \right)$ consists of exponentially many facet defining inequalities. Let us consider the following set:
\[
S^E = \left\{ (x, y, w) \in \mathbb{R}^{n + n + n}_+ : \sum_{i \in N} w_i \geq 1, w_i \leq l^{k_i} (x_i, y_i), k_i \in \bar{K}_i, i \in N, x \leq u \right\}.
\]

\begin{proposition}
The set $S^E$ is an extended formulation of $conv \left( S^U \right)$.
\end{proposition}

\begin{proof}
If $(x, y, w) \in S^E$, then clearly $(x, y) \in conv \left( S^U \right)$. Now let $(x, y) \in conv \left( S^U \right)$ and define $w_i = \min_{k_i} \{ l^{k_i} (x_i, y_i), k_i \in \bar{K}_i \}, i \in N$. Since the point $(x, y)$ is feasible for the set $conv \left( S^U \right), \sum_{i \in N} \min_{k_i} \{ l^{k_i} (x_i, y_i), k_i \in \bar{K}_i \} \geq 1$, and consequently $\sum_{i \in N} w_i \geq 1$. Thus $S^E$ is an extended formulation of $conv \left( S^U \right)$.
\end{proof}

Even though the description of $S^E$ consists of far fewer number of constraints than $conv \left( S^U \right)$, it has $\sum_{i \in N} u_i + n + 1$ linear inequalities in addition to the the bound constraints, which is pseudopolynomial in the input size because of its dependency on $u$.

This extended formulation can be solved as a linear program to optimize a linear function over $conv \left( S^U \right)$. When the components of $u$ are small (as in some cutting stock problems), this linear program can be solved fast.

\section{The Separation Problem}


We now describe a linear time separation algorithm to separate a given point $(\bar{x}, \bar{y})$ from $conv \left( S^U \right)$. Let $(\bar{x},\bar{y})$ be a point in $\mathbb{R}^n \times \mathbb{R}^n$. If $\bar{x} \nleq u$, then a bound constraint is sufficient to separate $(\bar{x},\bar{y})$. We thus consider the separation problem for the facet defining inequalities of $cl \left( conv \left( S^L \right) \right)$.

The facet defining inequalities of $cl \left( conv \left( S^L \right) \right)$ given by (\ref{FSL}) can be listed in a different way for easier understanding. Consider the following collection of columns.
\[
\begin{bmatrix}
    l^{1_1} (x_1,y_1)  & l^{1_2} (x_2,y_2) & l^{1_3} (x_3,y_3) & \dots & l^{1_n} (x_n,y_n) \\
    l^{2_1} (x_1,y_1)  & l^{2_2} (x_2,y_2) & l^{2_3} (x_3,y_3) & \dots & l^{2_n} (x_n,y_n) \\ \tag{M\textsubscript{U}} \label{Mu}
    \dots & \dots & \dots & \dots & \dots \\
    l^{(u_1 + 1)_1} (x_1,y_1)  & l^{(u_2 + 1)_2} (x_2,y_2) & l^{(u_3 + 1)_3} (x_3,y_3) & \dots & l^{(u_n + 1)_n} (x_n,y_n)
\end{bmatrix}
\]

Note that (\ref{Mu}) may have a different number of elements in each column
depending upon $u$, and thus it is not a matrix. The facet defining inequalities of $cl \left( conv \left( S^L \right) \right)$ can be constructed by adding $n$ terms from (\ref{Mu}), taking exactly one term from each column and constraining the sum to be at least one.

Let us revisit the example (\ref{E1}) in Section~\ref{Sec1}. As discussed in Section~\ref{Sec1}, the point $\left( 5, 1, 6, \frac{5}{6} \right)$ lies in $conv(S)$. But we see that this point is violated by the
inequality $\frac{5 y_1}{20} + \frac{6 y_2}{20} \geq 1$ which is of the form
(\ref{NF}). Note that the inequalities $\frac{5 y_1}{20} \geq 1$ and $\frac{6
y_2}{20} \geq 1$ are valid for $S_1$ and $S_2$, respectively and combining them in the way described above we obtain the inequality $\frac{5 y_1}{20} + \frac{6 y_2}{20} \geq 1$. Adding this inequality to $conv(S)$, we get the optimal solution $(5, 4, 6, 0)$ with optimal value 23.

If $(\bar{x},\bar{y}) \notin conv \left( S^U \right)$, then it must be violated by at least one inequality of the form (\ref{FSL}). In order to find such a violated inequality, we have to find one term from each column of (\ref{Mu}) so that the sum of these is less than 1.

\subsection{Efficient Separation for $conv \left( S^U \right)$}

In order to separate a point $(\bar{x}, \bar{y})$ from $conv \left( S^U \right)$, we find a minimum element from each column of (\ref{Mu}) at $(\bar{x}, \bar{y})$ and add them. Clearly, if the the sum is greater than or equal to 1, the point $(\bar{x},\bar{y})$ is feasible to $conv \left( S^U \right)$. Otherwise, adding the corresponding terms from each column and setting it to greater than or equal to 1, will give us a violated facet defining inequality.

Column $i$ of (\ref{Mu}) has $(u_i + 1)$ terms, $i \in N$. To solve the separation problem, we need to find the minimum value at $(\bar{x}, \bar{y})$ from each column. This step takes $O(u_i)$ time which is pseudo-polynomial in the size of input. We now present a linear time algorithm for the separation problem.

\begin{proposition}
There exists an efficient separation of the facet defining inequalities of $conv \left( S^U \right)$.
\end{proposition}

\begin{proof}
Since the bound constraints can be checked easily, let $(\bar{x}, \bar{y}) \in
\mathbb{R}^n_+ \times \mathbb{R}^n_+$ such that $\bar{x} \leq u$ be a given point. For each column of (\ref{Mu}), we want to find the term that gives the minimum evaluation at the point $(\bar{x}, \bar{y})$. Let
\[
\xi_i = \min \left\lbrace \frac{\bar{x}_i}{2w - 1} + \frac{\bar{y}_i w (w - 1)}{r (2w - 1)}, \frac{\bar{y}_i}{\bar{u}_i}, w = 1, \ldots, u_i \right\rbrace, \ \text{where } \bar{u}_i = \frac{r}{u_i}.
\]
Note that $\xi_i \geq 0$. To find $\xi_i$, we consider the following cases: \\
\textsc{Case 1:} If $\bar{y}_i = 0$, then clearly $\xi_i = 0$ at the last term, i.e., at $\frac{y_i}{\bar{u}_i}$ since $\frac{\bar{y}_i}{\bar{u}_i} = 0$. \\
\textsc{Case 2:} If $\bar{x}_i = 0$, then again $\xi_i = 0$ at $w = 1$. \\
\textsc{Case 3:} $\bar{x}_i > 0$ and $\bar{y}_i > 0$. Let us consider the following function:
\[
f(w) = \frac{\bar{x}_i}{2w - 1} + \frac{\bar{y}_i w (w - 1)}{r (2w - 1)}, w \geq 1.
\]

Our goal is to find a positive integer $q$ that minimizes $f(w)$ among all the integers in $[1, u_i]$. The function $f$ is continuously differentiable in the domain $w \geq 1$ with
\begin{align*}
f'(w) = - \frac{2 \bar{x}_i}{(2w - 1)^2} + \frac{\bar{y}_i}{r} \cdot \frac{2w^2 - 2w + 1}{(2w - 1)^2} \ \text{and } f''(w) & = \frac{2(4 \bar{x}_i r - \bar{y}_i)}{r(2w - 1)^3}.
\end{align*}
We have the following two subcases:

\textsc{Case 3.1:} When $4 \bar{x}_i r - \bar{y}_i > 0$, the function $f$ is strictly convex and has unique minimizer, say $\bar{w}_i$. Now $f'(\bar{w}_i) = 0$ occurs at
\begin{equation}
\bar{w}_i = \frac{1}{2} + \frac{\sqrt{\frac{4 \bar{x}_i r}{\bar{y}_i} - 1}}{2}. \label{BarWi}
\end{equation}

When $\bar{w}_i \leq 1$, the integer minimizer of $f$ is $q = 1$. When $u_i > \bar{w}_i > 1$, $q = \lceil \bar{w}_i \rceil$ or $\lfloor \bar{w}_i \rfloor$ whichever gives a lower $f(q)$ is the required $q$. Finally, $q = u_i$ when $\bar{w}_i > u_i$.

\textsc{Case 3.2:} When $4 \bar{x}_i r - \bar{y}_i \leq 0$, the function $f$ is concave for $w \geq 1$. Therefore, the minimum value will be attained at a boundary point, i.e., either at 1 or at $u_i$. Moreover, we see that
\begin{align*}
f'(w) & = - \frac{2 \bar{x}_i}{(2w - 1)^2} + \frac{\bar{y}_i}{r} \cdot \frac{2w^2 - 2w + 1}{(2w - 1)^2} \\
& = \frac{2 \bar{y}_i w(w - 1) + \bar{y}_i - 2 \bar{x}_i r}{r(2w - 1)^2} > 0.
\end{align*}

Thus, $f$ is strictly increasing function, and is minimized at $q = 1$. Now one more comparison is required to find the value of $\xi_i$. If $\frac{\bar{x}_i}{2q - 1} + \frac{\bar{y}_i q (q - 1)}{r(2q - 1)} \leq \frac{\bar{y}_i}{\bar{u}_i}$ then $\xi_i = \frac{\bar{x}_i}{2q - 1} + \frac{\bar{y}_i q (q - 1)}{r(2q - 1)}$, else $\xi_i = \frac{\bar{y}_i}{\bar{u}_i}$.

The term corresponding to any column of (\ref{Mu}) can be computed in $O(1)$ time, and since there are $n$ columns, a violated inequality can be found in $O(n)$ time.

If $\sum_{i=1}^n \xi_i \geq 1$, the point $(\bar{x}, \bar{y})$ is feasible to $conv \left( S^U \right)$. In Algorithm~\ref{Alg1} in Appendix~\ref{AppSx} we describe the separation algorithm in pseudocode.
\end{proof}

\begin{corollary}
The optimization problem having a linear objective function over $conv \left( S^U \right)$ can be solved in time polynomial in size of the input.
\end{corollary}

\begin{proof}
Since there is a polynomial time separation algorithm of the facet defining inequalities of $conv \left( S^U \right)$, the optimization of a linear function over $conv \left( S^U \right)$ can also be done in polynomial time (Gr{\"o}tschel et al. \cite{EquivSep}). We present an algorithm in Appendix~\ref{AppSx}. 
\end{proof}


\subsection{Efficient Separation for $conv(S)$}

The separation problem in the case of $conv(S)$ can also be solved in similar way with some modification. We use this algorithm to compare the effectiveness of our new cuts derived for $conv \left( S^U \right)$ in computational experiments.


\begin{proposition}
There exists an efficient separation of the facet defining inequalities of $conv(S)$.
\end{proposition}

\begin{proof}
Given a point  $(\bar{x}_i, \bar{y}_i) \in \mathbb{R}^n_+ \times \mathbb{R}^n_+$, let
\[
\xi_i = \min \left\lbrace \frac{\bar{x}_i}{2w - 1} + \frac{\bar{y}_i w (w - 1)}{r (2w - 1)}, w \in \mathbb{N} \right\rbrace.
\]

Note that $\xi \geq 0$ for $w \geq 1$. Our goal is to find a positive integer that minimizes $f(w)$. We consider the following cases.

\textsc{Case 1:} When $\bar{x}_i = 0$, then clearly $\xi_i = 0$ at $\hat{w}_i = 1$.

\textsc{Case 2:} When $\bar{y}_i = 0, \bar{x}_i \neq 0$, then  $\inf \left\lbrace \frac{\bar{x}_i}{2k - 1} + \frac{\bar{y}_i k (k - 1)}{r (2k - 1)}, k \in \mathbb{N} \right\rbrace = 0$, since $\frac{\bar{x}_i}{2k - 1} \rightarrow 0$ as $k \rightarrow \infty$. Therefore $\xi_i$ can be taken as 0 in this case.

\textsc{Case 3:} When $\bar{x}_i > 0, \bar{y}_i > 0$, the same logic used for $conv \left( S^U \right)$ can be deployed. Let $\hat{w}_i$ be the desired integer value. Then
\[
\hat{w}_i = 
\begin{cases}
 1, \ \text{when } 4 \bar{x}_i r - \bar{y}_i > 0 \ \text{and} \ \bar{w}_i \leq 1, \\
 
 \lceil \bar{w}_i \rceil,  \ \text{when} \ 4 \bar{x}_i r - \bar{y}_i > 0, \bar{w}_i > 1 \ \text{and } f \left( \lceil \bar{w}_i \rceil \right) \leq f \left( \lfloor \bar{w}_i \rfloor \right), \\
 \lfloor \bar{w}_i \rfloor,  \ \text{when} \ 4 \bar{x}_i r - \bar{y}_i > 0, \bar{w} > 1 \ \text{and } f \left( \lceil \bar{w}_i \rceil \right) \geq f \left( \lfloor \bar{w}_i \rfloor \right), \\
 
 1, \ \text{when } 4 \bar{x}_i r - \bar{y}_i \leq 0,
\end{cases}
\]
where $\bar{w}_i$ is defined by (\ref{BarWi}). If $\sum_{i=1}^n \xi_i \geq 1$, the point $(\bar{x}, \bar{y})$ is feasible to $conv(S)$. Otherwise, it is infeasible, and we have to find a violated facet defining inequality. We know the required value of $\hat{w}_i$ for \textsc{Case} 1 and 3. Let $t \in \mathbb{N}$ such that the following holds,
\begin{equation}
\sum_{i = 1}^n \xi_i + \sum_{i \in N : \bar{x}_i > 0, \bar{y}_i = 0} \frac{\bar{x}_i}{2t - 1} < 1. \label{Violation}
\end{equation}

Such a $t$ can always be found by the Archimedian property. A simple calculation shows that any integer greater than $\left\lfloor \frac{1 - \xi + v}{2(1 - \xi)} \right\rfloor$ where $\xi = \sum_{i = 1}^n \xi_i, v = \sum_{i \in N : \bar{x}_i > 0, \bar{y}_i = 0} \bar{x}_i$ is sufficient. Therefore, the following inequality is violated by the point $(\bar{x}, \bar{y})$:

\begin{align*}
\sum_{i \in N : \bar{x}_i = 0} x_i + \sum_{i \in N : \bar{x}_i > 0, \bar{y}_i > 0} \left[ \frac{x_i}{2\hat{w}_i - 1} + \frac{y_i \hat{w}_i (\hat{w}_i - 1)}{r(2 \hat{w}_i - 1)} \right] + \sum_{i \in N : \bar{x}_i > 0, \bar{y}_i = 0} \left[ \frac{x_i}{2t - 1} + \frac{y_i t (t - 1)}{r(2 t - 1)} \right] \geq 1,
\end{align*}
where $t \in \mathbb{N}$ such that $t \geq \left\lfloor \frac{1 - \xi + v}{2(1 - \xi)} \right\rfloor + 1$. In Algorithm~\ref{Alg2} in Appendix~\ref{AppS}, we provide the pseudocode of the separation algorithm. 
\end{proof}

Note that for any positive integer $t \geq \left\lfloor \frac{1 - \xi + v}{2(1
- \xi)} \right\rfloor + 1$, we get a violated inequality. From
(\ref{Violation}) we see that as $t$ increases, the violation also increases
and equals  $1 -\sum_{i = 1}^n \xi_i$ in the limiting case. Following this
argument, one may
conclude that the best inequality is the one with $t$ arbitrarily large.
However, this conclusion may not be correct because our measure of violation
is not normalized properly. Ideally we should find an inequality farthest from
the given point. Such a measure can be considered in future studies. 


\begin{corollary}
The optimization problem having a linear objective function over the set $S$ (or equivalently over $conv(S)$) can be solved in polynomial time.
\end{corollary}


We present a polynomial time algorithm to optimize a linear function over $conv(S)$ in Appendix~\ref{AppS}.


\section{Computational Results}

We now study the effectiveness of the cuts obtained for $S^U$ by doing computational experiments on cutting stock instances of the following form.
\begin{align*}
\min & \sum_{i=1}^n y_i \\
& \sum_{i \in N} x_{ij} y_i \geq d_j, \ j \in F,  \tag{CS} \label{CS} \\
& \sum_{j \in F} l_j x_{ij} \leq L, i \in N, \\
& x_{ij} \in \mathbb{Z}_+, y_i \in \mathbb{R}_+, \forall i \in N, j \in F,
\end{align*}
where the notation is the same as that in Section~\ref{Sec1}. These instances have stocks of one length $L$ from which $n$ different sizes of finals are to be cut. So, there are $n$ mixed-integer bilinear covering constraints modeling demand satisfaction. The upper bounds $x_{ij} \leq \left \lfloor \frac{L}{l_j} \right \rfloor = \nu_j (\text{say}), \forall i \in N, j \in F$ of the integral variables are implicit from the knapsack constraints present in the formulation. Here, our objective is to minimize the total number of stocks that are used. Since not more than $n$ finals are usually seen in solutions to (\ref{CS}), we assume $|N| = |F| = n$.

We have selected for our experiments ten instances used in Umetani et al.
\cite{Fiber} taken from applications in a chemical fiber company in Japan
(Fiber-xx-xxxx), six instances generated by CUTGEN (Gau and Wascher
\cite{CutGen}) (CutGen-xx-xx) and five randomly generated instances (Rand-xx).
These random instances were generated by fixing $L$ to 1030 and selecting
specifc problem size $n$ (denoted as `xx' in the name). The final lengths
$l_j$ were generated randomly between 75 and 600, and $d_j$ between 300 and
5000.  

We perform three sets of experiments. In all three we have used PuLP (Mitchell
et al. \cite{PuLP}) version 1.6.2 (installed in Python 2.7.12) to model the
linear programs and CBC (Forrest et al. \cite{Cbc}) solver to solve them. The
system we used to run our code has Linux (Ubuntu 16.04) operating system with
4x Intel(R) Core(TM) i5-3570 CPU@3.40 GHz processor and 8 GB of RAM. All
experiments were carried out on a single core.

In our first study we compare the bounds generated by our cuts for $conv \left( S^U \right)$ to those by Tawarmalani et al. \cite{OrtDis} for $conv(S)$. In both the cases we consider the facet defining inequalities of each mixed-integer bilinear covering constraint. Adding facet defining inequalities for each mixed-integer bilinear covering constraint together gives a polyhedral relaxation for the actual problem. For each instance, in both the cases, we start our iterations with the facet defining inequalities $\sum_{i=1}^n x_{ij} \geq 1$, for all $j \in F$, the bound constraints and the knapsack inequalities, i.e., we start our iterations by solving the following linear program.
\begin{align*}
\min & \sum_{i=1}^n y_i \\
\text{s.t. } & \sum_{i=1}^n x_{ij} \geq 1, \forall j \in F, \\
& 0 \leq x_{ij} \leq \nu_j, \forall i \in N, j \in F, \tag{LP-I} \label{LPI} \\
& \sum_{j \in F} l_j x_{ij} \leq L, i \in N, \\
& y \geq 0.
\end{align*}
Then, we add violated inequalities (if any) obtained from our separation
procedures and resolve (\ref{LPI}). This process is continued until we can not
find any more violated inequalities, or the number of LPs solved exceeds a
predefined limit of 800, or the total time used exceeds two hours. If we can
not find any more violated inequalities, then the solution of the current LP
lies in the convex hull of the set $S^U$ associated with each of the bilinear
constraints. This solution may not be feasible to the original
problem~(\ref{CS}).

We run the above experiment in two different settings using facet defining inequalities derived (i) for $conv \left( S^U \right)$ and (ii) for $conv(S)$. We consider the sets $S^U$ and $S$ by looking at each bilinear covering constraint separately and add one most violated cut for each such constraint using Algorithm~\ref{Alg1} and Algorithm~\ref{Alg2} respectively. So, we add at most $n$ cuts in every iteration (LP solve) which are not deleted in further iterations. This means, at iteration $k$, we solve an LP relaxation of the  instance with at most $k|F|$ number of linear inequalities in addition to those in (\ref{LPI}).

\begin{figure}
\centering
\begin{minipage}{.5\textwidth}
  \centering
	\includegraphics[width=\textwidth]{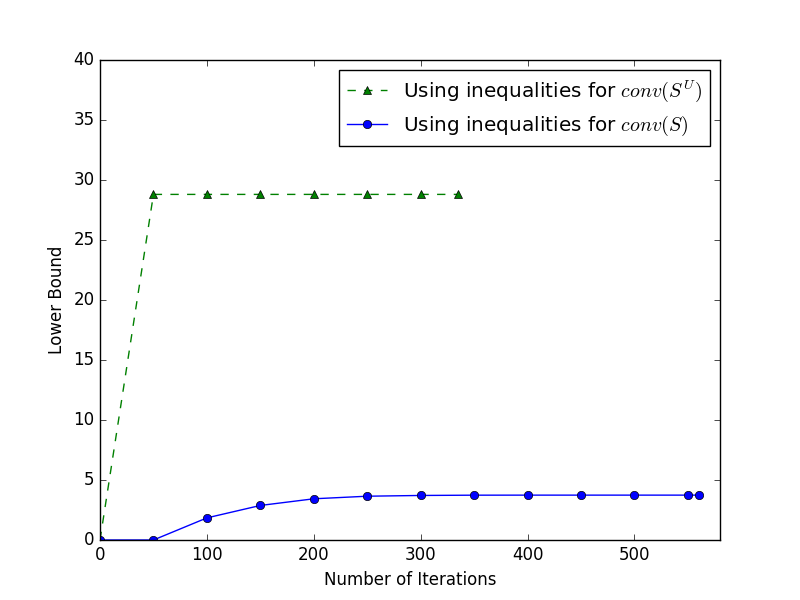}
	\caption{Bound comparisons for Fiber-15-5180} \label{ComGraphFiber}
\end{minipage}%
\begin{minipage}{.5\textwidth}
  \centering
	\includegraphics[width=\textwidth]{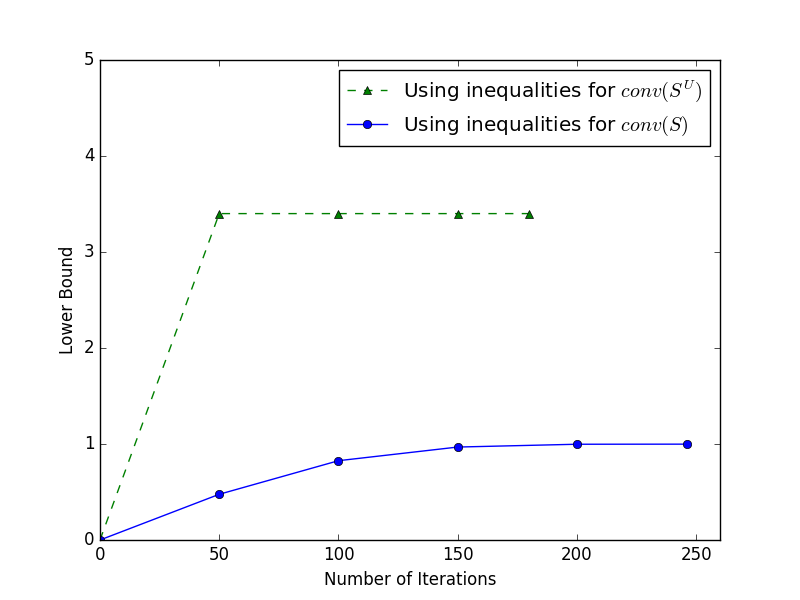}
	\caption{Bound comparisons for CutGen-01-25} \label{ComGraphCutGen}
\end{minipage} \\
\includegraphics[scale=0.32]{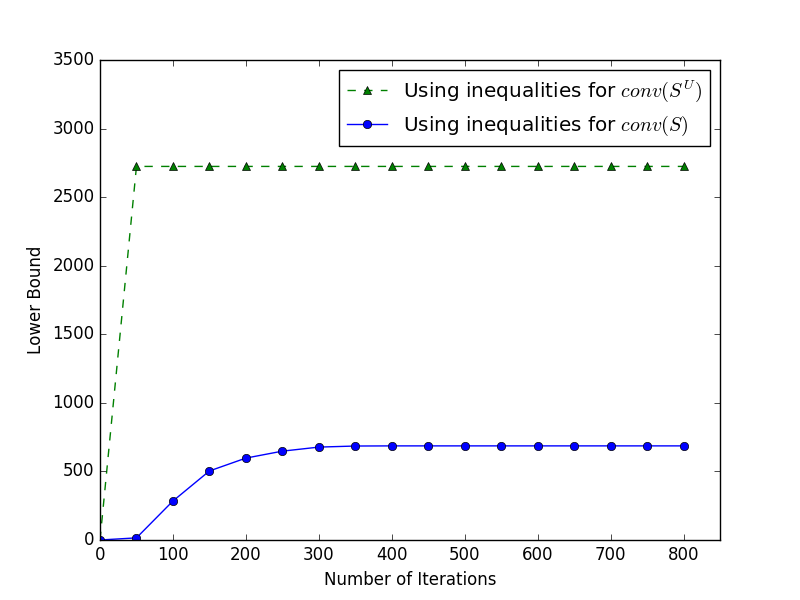}
\caption{Bound comparisons for Rand16} \label{ComGraphRand}
\end{figure}

Table~\ref{Table1} compares the effects of new cuts for $S^U$ to the cuts
derived for $S$. We observe that cuts for $conv \left( S^U \right)$ improve
the lower bounds with fewer cuts and in lesser time as compared to cuts for
 $conv(S)$. In the Figures~\ref{ComGraphFiber},
\ref{ComGraphCutGen} and \ref{ComGraphRand}, we present iteration-wise bound
comparisons for three instances Fiber-15-5180, CutGen-01-25 and Rand16
respectively. We see that the cuts for $conv \left( S^U \right)$ improve
bounds faster than those for $conv(S)$.

We also study the time taken to solve the extended formulations (both LP and
MILP) of Section~\ref{SecExtFor}. While the LP defines the set $conv \left( S^U \right)$, the MILP is an even tighter relaxation of cutting stock problem. Recall that the extended formulation has $2
n^2 + n$ variables, that means $n^2$ more variables than the original
formulation. Table~\ref{Table2} lists the bounds and time taken to solve the
two relaxations. We set computational time limit
to two hours. We write ``7200*" for the instances where this time limit is
reached, and for such instances we report the relative gap $\frac{ub - lb}{lb}$ of the MILP. The $lb$ of the MILP is a lower bound for the optimal value of (\ref{CS}).  We
also compute an upper bound to optimal solution of~(\ref{CS}) obtained by
fixing the  variable $x$ to the MILP
solution in (\ref{CS}) and solving a linear program in $y$ only. This bound is
reported in the last column (``UB'') of Table~\ref{Table2}.

We see that the extended
formulation LP takes much less time compared to the cutting plane algorithm
using cuts for $conv \left( S^U \right)$, and even the MILP is often faster
than the cut
based iterative LP approach. This observation suggests that extended formulation is quite good
for these instances when the implied bounds $\nu$ on $x$ are small. The extended MILP for randomly generated instances seems to be unusually difficult for the solver. The bounds
given by the LP and MILP of the extended formulation are the same for all
instances except for Rand10. We do not have an explanation of this phenomenon
currently.

\begin{table}
\centering
\begin{small}
\caption{Comparison of iterations taken to optimize over the convex hull and the lower bounds obtained. (Here ``Iter" means number of LP iterations, and ``LB" means Lower Bound obtained after termination, ``Cuts" column indicates the number of cuts added, Time is in seconds). A * mark indicates time or iteration limit is reached} \label{Table1}
\begin{tabular}{|c|c|c|c|c|c|c|c|c|c|}
\hline
\multirow{2}{*}{Instance} & \multirow{2}{*}{$n$} & \multicolumn{4}{|c|}{Using inequalities for $conv \left( S^U \right)$} & \multicolumn{4}{|c|}{Using inequalities for $conv(S)$} \\
\cline{3-10}
& & Iter & Cuts & LB & Time & Iter & Cuts & LB & Time \\
\hline
Fiber10-5180 & 10 & 156 & 1212  & 27.00    & 18.23 & 226 & 1917    & 6.88   & 55.09 \\
\hline
Fiber10-9080 & 10 & 215 & 1663  & 15.00    & 23.86 & 223 & 2045    & 3.85   & 46.15 \\
\hline
Fiber11-5180 & 11 & 118 & 939   & 26.00    & 9.22  & 288 & 2673    & 6.10   & 89.22 \\
\hline
Fiber11-9080 & 11 & 463 & 3151  & 14.44    & 137.83 & 335 & 2946   & 3.40   & 162.3 \\
\hline
Fiber14-5180 & 14 & 147 & 1343  & 22.00    & 17.27 & 473 & 5417    & 3.34   & 547.81 \\
\hline
Fiber14-9080 & 14 & 136 & 1522  & 11.00    & 19.22 & 476 & 6211    & 1.90   & 658.27 \\
\hline
Fiber15-5180 & 15 & 335 & 2350  & 28.80    & 98.96 & 560 & 7219    & 3.74   & 1412.46 \\
\hline
Fiber15-9080 & 15 & 623 & 2861  & 16.00    & 317.27 & 800* & 8890  & 2.09   & 1881.71 \\
\hline
Fiber16-5180 & 16 & 800* & 2393 & 27.20    & 282.51 & 756 & 9763   & 5.17   & 3086.56 \\
\hline
Fiber16-9080 & 16 & 223 & 2566  & 15.11    & 63.9 & 723 & 10330    & 2.93   & 2692.61 \\ 
\hline
CutGen01-01 & 10 & 211 & 1359   & 2.43	   & 30.79 & 252 & 2244    & 1.24   & 95.33 \\
\hline
CutGen01-02 & 10 & 235 & 1888   & 2.57     & 43.86 & 270 & 2443    & 0.97   & 103.32 \\
\hline
CutGen01-25 & 10 & 180 & 1238   & 3.40     & 21.95 & 246 & 2157    & 0.99   & 77.18 \\
\hline
CutGen01-100 & 10 & 194 & 1276  & 3.80     & 18.82 & 244 & 2131    & 1.25   & 52.82 \\
\hline
CutGen02-40 & 10 & 186 & 1292   & 26.00    & 22.48 & 262 & 2272    & 10.41  & 	92.60 \\
\hline
CutGen02-60 & 10 & 186 & 1322   & 33.80    & 24.05 & 275 & 2480    & 10.10  & 76.32 \\
\hline
Rand10      & 10 & 64 & 557     & 1520.50  & 3.06 & 185 & 1601     & 697.22 & 28.86 \\
\hline
Rand15      & 15 & 114 & 1272   & 2122.00  & 18.15 & 792 & 8485    & 576.69 & 3786.64 \\
\hline
Rand16      & 16 & 800* & 2293  & 2724.00  & 257.04 & 800* & 8780  & 686.27 & 4611.77 \\
\hline
Rand20      & 20 & 173 & 2698   & 2250.00  & 66.91 & 725 & 13131   & 631.02 & 7200* \\
\hline
Rand25      & 25 & 800* & 9835  & 1437.00  & 3382.22 & 686 & 17175 & 517.15 & 7200* \\
\hline
\end{tabular}
\end{small}
\end{table}

\begin{table}
\centering
\begin{small}
\caption{Comparison for Extended Formulation. (``Const." column contains total number of constraints in the extended formulation)} \label{Table2}
\begin{tabular}{|c|c|c|c|c|c|c|c|}
\hline
\multirow{2}{*}{Instance} & \multirow{2}{*}{Const.} &
\multicolumn{2}{|c|}{Extended LP} & \multicolumn{3}{|c|}{Extended MILP} & \multirow{2}{*}{UB} \\
\cline{3-7}
& & LB & Time & LB(Rel. Gap) & Time & Nodes &  \\
\hline
Fiber10-5180 & 610 & 27.00	& 0.17 & 27.00 & 1.64 & 41 & 135.00 \\
\hline
Fiber10-9080 & 1010 & 15.00	& 0.25 & 15.00 & 2.99 & 61 & 68.06 \\
\hline
Fiber11-5180 & 792 & 26.00 & 0.22 & 26.00 & 6.21 & 91 & 75.28 \\
\hline
Fiber11-9080 & 1331 & 14.44 & 0.37 & 14.44 & 3.55 & 4 & 45.20 \\
\hline
Fiber14-5180 & 1274 & 22.00 & 0.34 & 22.00 & 4.32 & 176 & 66.00 \\
\hline
Fiber14-9080 & 2114 & 11.00 & 0.59 & 11.00 & 1.68 & 1 & 40.86 \\
\hline
Fiber15-5180 & 1470 & 28.80 & 0.57	 & 28.80 & 5.80 & 91 & 88.00 \\
\hline
Fiber15-9080 & 2430 & 16.00 & 0.86 & 16.00 & 8.54 & 91 & 35.02 \\
\hline
Fiber16-5180 & 1648 & 27.20 & 0.69	 & 27.20 & 9.09 & 615 & 136.00 \\
\hline
Fiber16-9080 & 2800 & 15.11 & 0.74	 & 15.11 & 7.38 & 171 & 82.50 \\ 
\hline
CutGen01-01 & 1740 & 2.43 & 0.68 & 2.43 & 18.78 & 635 & 13.25 \\
\hline
CutGen01-02 & 1300 & 2.57 & 0.55 & 2.57 & 9.76 & 208 & 10.77 \\
\hline
CutGen01-25 & 1550 & 3.40 & 0.44 & 3.40 & 6.19 & 91 & 17.00 \\
\hline
CutGen01-100 & 1190 & 3.80 & 0.32 & 3.80 & 7.62 & 146 & 19.00 \\
\hline
CutGen02-40	 & 2170 & 26.00 & 0.67 & 26.00 & 13.56 & 116 & 149.00 \\
\hline
CutGen02-60 & 1780 & 33.80 & 0.53 & 33.80 & 6.58 & 80 & 117.36 \\
\hline
Rand10 & 320 & 1520.50 & 0.11 & 1557.87(1.46) & 7200* & 7562839 & 7407.00 \\
\hline
Rand15 & 720 & 2122.00 & 0.27 & 2122.00(1.89) & 7200* & 3818886 & 9688.00 \\
\hline
Rand16 & 832 & 2724.00 & 0.30 & 2724.00(1.58) & 7200* & 2392670 & 12899.00 \\
\hline
Rand20 & 1440 & 2250.00 & 0.53 & 2250.00(2.59) & 7200* & 1758597 & 13990.00 \\
\hline
Rand25 & 3750 & 1437.00 & 1.67 & 1437.00(3.67)	& 7200* & 681074 & 15041.25 \\
\hline
\end{tabular}
\end{small}
\end{table}

Lastly, we consider an exact MILP formulation of (\ref{CS}). Let $w_{ijh} = 1$ if $x_{ij} = h$ and $w_{ijh} = 0$ otherwise for $i \in N, j \in F$ and $h \in \left\{0, 1, \ldots, \nu_j \right\}$. Replacing the terms $w_{ijh} y_i$ with $z_{ijh}$ and using the linear inequalities to model the products $z_{ijh} y_i$ we have the following formulation:
\begin{align*}
\min & \sum_{i=1}^n y_i \\
& \sum_{i \in N} \sum_{h = 0}^{\nu_j} h z_{ijh} \geq d_j, \ j \in F,  \\
& \sum_{j \in F} l_j x_{ij} \leq L, i \in N, \\
& \sum_{h = 0}^{\nu_j} w_{ijh} = 1, i \in N, j \in F, \tag{R\textsubscript{CS}} \label{RCS} \\
& \sum_{h = 0}^{\nu_j} h w_{ijh} = x_{ij}, i \in N, j \in F,  \\
& z_{ijh} \geq y_i + B w_{ijh}, i \in N, j \in F, h \in  \left\{0, 1, \ldots, \nu_j \right\}, \\
& z_{ijh} \leq y_i, i \in N, j \in F, h \in \left\{0, 1, \ldots, \nu_j \right\}, \\
& z_{ijh} \leq B w_{ijh}, i \in N, j \in F, h \in \left\{0, 1, \ldots, \nu_j \right\}, \\
& w_{ijh} \in \{0, 1\}, z_{ijh} \in \mathbb{R}_+, x_{ij} \in \mathbb{Z}_+, y_i \in \mathbb{R}_+, \forall i \in N, j \in F,  \in \left\{0, 1, \ldots, \nu_j \right\}.
\end{align*}

The formulation (\ref{RCS}) is an exact reformulation of (\ref{CS}) because
$w_{ijh}$ are binary. Here $B$ is an upper bound for the variables $y$. For
our experiment, we used the UB value reported in Table~\ref{Table2} for $B$.

In Table~\ref{Table3}, we list both lower and upper bounds to the objective
value of (\ref{CS}) by solving the MILP reformulation (\ref{RCS}). We also
report the lower bounds obtained from cuts for $conv(S^U)$ and $conv(S)$ from
the earlier tables for comparison. The time limit was again set to two hours.
We observe that the solver reached the time limit for all
instances. Further the lower bound at the root relaxation was the same as the
lower bound after two hours for all instances. From the table we see that the
lower bounds generated by solving (\ref{RCS}) are smaller than the lower
bounds generated by the facets of $conv \left( S^U \right)$ in all instances
except Rand15, Rand16 and Rand20 for which they are equal. On the other hand,
the bounds from $conv(S)$ are sometimes weaker than that from the
MILP~(\ref{RCS}).

\begin{table}
\centering
\begin{small}
\caption{Bounds generated by Binary MILP (\ref{RCS}) after two hours of computational time.} \label{Table3}
\begin{tabular}{|c|c|c|c|c|p{12ex}|p{10ex}|}
\hline
Instance & LB & UB & Rel. Gap & Nodes & LB ($conv(S^U)$) & LB ($conv(S)$) \\
\hline
Fiber10-5180 & 9.00 & 69.80 & 6.76 & 754875     & 27.00    & 6.88              \\
\hline                                                              
Fiber10-9080 & 3.00 & 39.19 & 12.06 & 232134    & 15.00    & 3.85              \\
\hline                                                              
Fiber11-5180 & 8.67 & 67.15 & 6.75 & 651330     & 26.00    & 6.10              \\
\hline                                                              
Fiber11-9080 & 2.89 & 39.79 & 12.77 & 216931    & 14.44    & 3.40              \\
\hline                                                              
Fiber14-5180 & 11.00 & 48.83 & 3.44 & 311882    & 22.00    & 3.34              \\
\hline                                                              
Fiber14-9080 & 3.14 & 32.44 & 9.32 & 91340      & 11.00    & 1.90              \\
\hline                                                              
Fiber15-5180 & 9.60 & 60.73 & 5.33 & 220926     & 28.80    & 3.74              \\
\hline                                                              
Fiber15-9080 & 3.20 & 34.51 & 9.78 & 51267      & 16.00    & 2.09              \\
\hline                                                              
Fiber16-5180 & 9.50 & 85.97 & 8.05 & 347264     & 27.20    & 5.17              \\
\hline                                                              
Fiber16-9080 & 3.39 & 66.67 & 18.65 & 78161     & 15.11    & 2.93              \\ 
\hline                                                              
CutGen01-01 & 0.62 & 13.56 & 20.90 & 215080     & 2.43	   & 1.24              \\
\hline                                                              
CutGen01-02 & 0.64 & 10.62 & 15.52 & 138215     & 2.57     & 0.97              \\
\hline                                                              
CutGen01-25 & 1.13 & 10.69 & 8.43 & 119516      & 3.40     & 0.99              \\
\hline                                                              
CutGen01-100 & 1.27 & 13.40 &  9.58 & 210161    & 3.80     & 1.25              \\
\hline                                                              
CutGen02-40 & 8.67 & 117.97 & 12.61 & 114179    & 26.00    & 10.41             \\
\hline                                                              
CutGen02-60 & 11.27 & 114.66 & 9.18 & 136551    & 33.80    & 10.10             \\
\hline                                                              
Rand10 & 1013.67 & 7101.50 & 6.01 & 1501923     & 1520.50  & 697.22            \\
\hline                                                              
Rand15 & 2122.00 & 8834.67 & 3.16 & 554079      & 2122.00  & 576.69            \\
\hline                                                              
Rand16 & 2724.00 & 10876.00 & 2.99 & 1310800    & 2724.00  & 686.27            \\
\hline                                                              
Rand20 & 2250.00 & 13608.83 & 5.05 & 183486     & 2250.00  & 631.02            \\
\hline                                                              
Rand25 & 958.00 & 14200.16 & 13.82 & 103304     & 1437.00  & 517.15            \\
\hline
\end{tabular}
\end{small}
\end{table}


\section{Concluding Remarks and Future Work}

When bounds on integer variables in a bilinear covering set are finite, we
are able to obtain the polyhedral description of 
the convex hull. Even though one can not directly apply the orthogonal
disjunctive
procedure here, we are still able to compute the convex hull by
first creating a suitable relaxation and then applying the procedure. It would be
interesting to see if similar procedures can be applied to other restrictions
of the set as well.  Our examples and experiments show that the new facet
defining inequalities of $conv \left( S^U \right)$ improve the bounds as compared to the case
when bounds are not considered. The extended formulation for many cutting stock problem can be solved fast, even the MILP can be solved fast.

The procedure of finding facet defining
inequalities to separate a given point from the convex hull is fast.
Our results can be applied in a straight-forward manner to
the following set also:
\[
S^\delta = \left\{ (x,y) \in \mathbb{Z}^n_+ \times \mathbb{R}^n_+ : \sum_{i=1}^n \delta_i x_i y_i \geq r, x \leq u \right\},
\]
where $r > 0, u \in \mathbb{N}$ and $\delta_i > 0$ for all $i \in N$. Using
our analysis, we can show that the facet defining inequalities of $conv \left( S^\delta
\right)$ also can be separated in $O(n)$ time.

In this work we did not consider the knapsack constraint of (\ref{CS}) in our set. Including
it and also considering multiple bilinear constraints together can be taken up
as future work. 

The cut generated by our criterion of `maximum violation' without any
normalization  may not be the cut that improves the lower bound the most, or the cut that is farthest from the infeasible point.
Consider the following example:
\begin{align*}
\min \ x_1 + y_1 + x_2 & + y_2 \\
\text{s.t. }  x_1 y_1 + x_2 y_2 & \geq 20, \\
 x_i & \leq 10, \quad i = 1, 2, \\
 x_i \in \mathbb{Z}_+, y_i & \geq 0,  \quad i = 1, 2.
\end{align*}

The point $(x_1, y_1, x_2, y_2) = (5,4,0,0)$ is a global minimizer with
optimal value 9. At the first iteration, the
LP solution is $(1, 0, 0, 0)$ with objective value 1. The best cut generated by
Algorithm~\ref{Alg1} to cut this point off is $\frac{y_1}{2} + \frac{y_2}{2}
\geq 1$. After adding this, the solution is $(1,2,0,0)$ with objective value
3. But, if we instead add the
facet defining inequality $\frac{x_1}{5} + \frac{3 y_1}{50} + \frac{x_2}{5} +
\frac{3 y_2}{50} \geq 1$, we get a better solution $(5,0,0,0)$ with objective
value 5. Also, the distance of the latter from the point $(1, 0, 0, 0)$ is
nearly 2.7 as compared to 1.41 for the former.
Finding the cut that improves the bound the most or that is
farthest from the given point is another interesting question.




%


%
%
%

\section*{Acknowledgments} We are grateful to three anonymous reviewers for giving insightful comments and suggestions, specially the idea of the extended formulation and simplification of proofs. These suggestions have improved this paper significantly. We also thank Vishnu Narayanan for insightful discussions we had with him during this work.


\begin{appendices}

\section*{Appendices}

\section{Optimization Over $S$ and Separation on $conv(S)$} \label{AppS}

Pseudocode for separation of facets of $conv \left(S \right)$ is provided
in Algorithm~\ref{Alg2}. Now we consider the problem of  
minimizing a linear function $c^T x + d^T y$ over $S$ (or equivalently over $conv(S)$).

\begin{algorithm}
\caption{Separation of the facet defining inequalities of $conv(S)$} \label{Alg2}
\begin{algorithmic}[1]
\State Input : A point $(\bar{x}, \bar{y}) \in \mathbb{R}^n_+ \times \mathbb{R}^n_+$
\State Output : Decide whether $(\bar{x}, \bar{y}) \in conv(S)$, and if not then provide a facet defining inequality that cuts off $(\bar{x}, \bar{y})$

\For{$i = 1, \ldots, n$}
	\If{$\bar{x}_i = 0$} \State $\hat{w}_i = 1, \xi_i = 0$
	\ElsIf{$\bar{x}_i \bar{y}_i > 0$}
		\If{$4 \bar{x}_i r > \bar{y}_i$}
			\If{$\frac{1}{2} + \frac{\sqrt{\frac{4 \bar{x}_i r}{\bar{y}_i} - 1}}{2} > 1$}
				\State $p = \left\lfloor \frac{1}{2} + \frac{\sqrt{\frac{4 \bar{x}_i r}{\bar{y}_i} - 1}}{2} \right\rfloor$
				\If{$\frac{\bar{x}_i}{2p - 1} + \frac{\bar{y}_i p (p - 1)}{r(2p - 1)} \leq \frac{\bar{x}_i}{2(p+1) - 1} + \frac{\bar{y}_i p (p + 1)}{r(2(p+1) - 1)}$}
					\State $\hat{w}_i = p$					
				\Else \State $\hat{w}_i = p + 1$
				\EndIf
			\Else \State $\hat{w}_i = 1$
			\EndIf
		\Else \State $\hat{w}_i = 1$
		\EndIf
		\State $\xi_i = \frac{\bar{x}_i}{2 \hat{w}_i - 1} + \frac{\bar{y}_i \hat{w}_i (\hat{w}_i - 1)}{r(2\hat{w}_i - 1)}$
	\Else \State $\xi_i = 0$	
	\EndIf
\EndFor
\State $\xi = \sum_{i \in N} \xi_i$
\If{$\xi \geq 1$}
	\State The point $(\bar{x}, \bar{y})$ is feasible to $conv(S)$.
\Else
	\State $v = \sum_{i \in N : \bar{y}_i = 0} \bar{x}_i$
	\State $t = \left\lfloor \frac{1 - \xi + v}{2(1 - \xi)} \right\rfloor + \gamma$, where $\gamma$ can be taken as any positive integer.
	\For{$i = 1, \ldots, n$}
		\If{$\bar{x}_i > 0$ and $\bar{y}_i = 0$}
			\State $\hat{w}_i = t$
		\EndIf
	\EndFor
	\State The inequality $\sum_{i = 1}^n \frac{x_i}{2 \hat{w}_i - 1} + \frac{y_i \hat{w}_i (\hat{w}_i - 1)}{r(2\hat{w}_i - 1)} \geq 1$ cuts off the point $(\bar{x}, \bar{y})$.
\EndIf
\end{algorithmic}
\end{algorithm}

If one of the components of $c$ or $d$ is negative, then the problem is unbounded. Suppose, $c \geq 0, d \geq 0$ and one of the component of the vector $c, c_t$ (say) is zero. If $d_t = 0, \min_{(x,y) \in S} c^T x + d^T y = 0$ and if $d_t > 0, \inf_{(x,y) \in S} \ c^T x + d^T y = 0$. This is because, in either case we can choose $y_t$ arbitrary small such that $x_t y_t = r$ and all other components are zero. Now, let $c \geq 0, d = 0$. Let $c_t \leq c_j, \forall j \in N$. Then $\mathcal{L}(t, 1, r)$ is an optimal solution with optimal value $c_t$.

 The only remaining case is when $c > 0, d \geq 0, d \neq 0$. We consider it next.

\begin{proposition} \label{2DPol}
Consider the orthogonal disjunctive subset $S_i$ of the set $S$. Then we can solve the optimization problem $\min_{(x,y) \in S_i} c_i x_i + d_i y_i$ in polynomial time.
\end{proposition}

\begin{proof}
From the definition, each $(x,y) \in S_i$ is of the form $\mathcal{L}(i, x_i, y_i), x_i \in \mathbb{N}$. If $c_i \geq 0, d_i = 0$, then $\mathcal{L}(i, 1, r)$ is an optimal solution with optimal value $c_i$.

Now, we only have to consider $c_i > 0, d_i > 0$. Let $\mathcal{L}(i, x^*_i, y^*_i)$ be an extreme point optimal solution of $conv \left( S_i \right)$. Clearly, this point should lie on the surface $x_i y_i = r$. Since the continuous relaxation of the set $S_i$ is a strictly convex set, the optimal solution $\mathcal{L}(i, \bar{x}_i, \bar{y}_i)$ (say) over the continuous relaxation is unique, and we have,
\[
\bar{x}_i = \sqrt{\frac{r d_i}{c_i}}, \ \text{and} \ \bar{y}_i = \frac{r}{\bar{x}_i}.
\]

If $\sqrt{\frac{r d_i}{c_i}}$ is an integer, then $x^*_i =  \bar{x}_i, y^*_i = \bar{y}_i$. If not, then from the geometry, it is clear that at the optimal solution either $x^*_i = \left\lceil \sqrt{\frac{r d_i}{c_i}} \right\rceil$ or $x^*_i = \left\lfloor \sqrt{\frac{r d_i}{c_i}} \right\rfloor$ whichever minimizes the objective function and is nonzero.

So, to find an optimal solution, we just have to check the signs of the objective coefficient and compute the value of $\sqrt{\frac{r d_i}{c_i}}$. This can be done in constant time. 
\end{proof}

Now we consider the set $S$. If an optimal solution exists, there must be an extreme point optimal solution of $conv(S)$ that is optimal in $S$. Now by Theorem~\ref{ExtThm}, there must be an optimal solution that is an extreme point of $conv(S_i)$ for some $i \in N$. We can solve the $n$ problems $\min_{\mathcal{L}(i, x_i, y_i) \in S_i} c_i x_i + d_i y_i$ for $i \in N$ and pick the minimum of the $n$ objective values, we will get the optimal value and corresponding optimal solution. Since each subproblem takes constant time to solve, we can solve the whole problem in linear time.



\section{Optimization over $S^U$ and Separation on $conv \left(S^U \right)$} \label{AppSx}

Pseudocode for separation of facets of $conv \left(S^U \right)$ is provided
in Algorithm~\ref{Alg1}. 
Now we consider the following problem:
\[
\zeta = \min_{(x,y) \in S^U} c^T x + d^Ty \tag{$P$} \label{P}
\]

\begin{algorithm}
\caption{Separation of the facet defining inequalities of $conv \left( S^U \right)$} \label{Alg1}
\begin{algorithmic}[1]
\State Input : A point $(\bar{x}, \bar{y}) \in \mathbb{R}^n_+ \times \mathbb{R}^n_+, x \leq u$
\State Output : Decide whether $(\bar{x}, \bar{y}) \in conv \left( S^U \right)$, and if not then provide a facet defining inequality of $conv \left( S^U \right)$ that cuts off $(\bar{x}, \bar{y})$

\For{$i = 1, \ldots, n$}
	\If{$\bar{y}_i = 0$} \State $\hat{w}_i = u_i + 1$
	\ElsIf{$\bar{x}_i = 0$} \State $\hat{w}_i = 1$
	\Else
		\State $q = 0$
		\If{$4 \bar{x}_i r > \bar{y}_i$}
			\If{$\frac{1}{2} + \frac{\sqrt{\frac{4 \bar{x}_i r}{\bar{y}_i} - 1}}{2} > 1$}
				\If{$\frac{1}{2} + \frac{\sqrt{\frac{4 \bar{x}_i r}{\bar{y}_i} - 1}}{2} < u_i$}
					\State $p = \left\lfloor \frac{1}{2} + \frac{\sqrt{\frac{4 \bar{x}_i r}{\bar{y}_i} - 1}}{2} \right\rfloor$
						\If{$\frac{\bar{x}_i}{2p - 1} + \frac{\bar{y}_i p (p - 1)}{r(2p - 1)} \leq \frac{\bar{x}_i}{2(p+1) - 1} + \frac{\bar{y}_i p (p + 1)}{r(2(p+1) - 1)}$}
							\State $q = p$					
						\Else \State $q = p + 1$
						\EndIf
				\Else \State $q = u_i$
				\EndIf
			\Else \State $q = 1$
			\EndIf
		\Else \State $q = 1$
		\EndIf
		\If{$\frac{\bar{x}_i}{2q - 1} + \frac{\bar{y}_i q (q - 1)}{r(2q - 1)} \leq \frac{\bar{y}_i}{\bar{u}_i}$}
			\State $\hat{w}_i = q$
		\Else \State $\hat{w}_i = u_i + 1$
		\EndIf
	\EndIf
\EndFor
\State $R = \sum_{i \in N : \hat{w}_i \leq u_i} \frac{\bar{x}_i}{2 \hat{w}_i - 1} + \frac{\bar{y}_i \hat{w}_i (\hat{w}_i - 1)}{r(2\hat{w}_i - 1)} + \sum_{i \in N : \hat{w}_i = u_i + 1} \frac{\bar{y}_i}{\bar{u}_i}$
\If{$R \geq 1$}
	\State The point $(\bar{x}, \bar{y})$ is feasible to $conv \left( S^U \right)$.
\Else
	\State The inequality $\sum_{i \in N : \hat{w}_i \leq u_i} \frac{x_i}{2 \hat{w}_i - 1} + \frac{y_i \hat{w}_i (\hat{w}_i - 1)}{r(2\hat{w}_i - 1)} + \sum_{i \in N : \hat{w}_i = u_i + 1} \frac{y_i}{\bar{u}_i} \geq 1$ separates $(\bar{x}, \bar{y})$.
\EndIf
\end{algorithmic}
\end{algorithm}

This problem is equivalent to minimizing $c^T x + d^Ty$ over $conv \left( S^U \right)$ which is polyhedral and whose extreme points are known. If an optimal solution exists, we will find an extreme point optimal solution to $conv \left( S^U \right)$.



When $d_t < 0$ for some $t \in N$, the problem is unbounded. Otherwise, if $c \leq 0, d = 0$, then clearly $\left( u_1, \frac{r}{u_1}, u_2, 0, \ldots, u_n, 0 \right)$ is an extreme point optimal solution.

Now the remaining case is $d \geq 0$. We first partition the set of extreme points of $conv \left( S^U \right)$ and optimize over those partitions. Let us define the following set for each $i \in N$.
\[
E_i = \left\{ (x,y) \in \mathbb{R}^{2n}_+ : x_i \in \{1, \ldots, u_i \}, y_i = \frac{r}{x_i}, x_j \in \{ 0, u_j \}, y_j = 0, \forall j \in N, j \neq i \right\}.
\]

From the discussion in Section~\ref{ExtDes}, all the points in $E_i$ are extreme points of $conv \left( S^U \right)$ and $E = \bigcup_{i \in N} E_i$ is the set of all extreme points of $conv \left( S^U \right)$. We minimize $c^T x + d^Ty$ over each set $E_i, i \in N$ and pick the minimum. Now our goal is to solve the following problem.
\[
\zeta_i = \min_{(x,y) \in E_i} c^T x + d^Ty \tag{$P_i$} \label{Pi}
\]

Clearly $\zeta = \min \{ \zeta_i : i \in N \}$. Note that only the $i^{th}$ component of the variable $y$ of each point in $E_i$ is non-zero and rest are all zero. Therefore, the objective function of the above problem (\ref{Pi}) reduces to $c_i x_i + d_i  y_i + \sum_{j \in N, j \neq i} c_j x_j$. For any point $(x,y) \in E_i$, the choices of the  components $x_j \in \{0, u_j\}, j \in N, j \neq i$ are independent of the choice of $x_i \in \{1, \ldots, u_i \}$. Let $(\bar{x}, \bar{y})^i \in E_i$ be an optimal solution of (\ref{Pi}). Then we must have $\bar{y}^i_i = \frac{r}{\bar{x}^i_i}, \bar{y}^i_j = 0, \forall j \in N, j \neq i$. Let us consider the following choices of $x$ components of $(\bar{x}, \bar{y})^i$.
\begin{align*}
& \bar{x}^i_i \in \{1, \ldots, u_i \} \ \text{such that } (\bar{x}^i_i, \bar{y}^i_i) \ \text{minimzes }  c_i x_i + d_i y_i, \\
& \bar{x}^i_j = \begin{cases} 0, \ \text{if } c_j > 0, \\ u_j, \ \text{if } c_j \leq 0,   \end{cases} \forall j \in N, j \neq i.
\end{align*}

It can be seen clearly that such above choice of the components of $(\bar{x}, \bar{y})^i$ minimizes the objective function. Now to find the value of $\bar{x}^i_i \in \{1, \ldots, u_i \}$, we consider the following cases.

\textsc{Case 1:} When $c_i \leq 0$, then $\bar{x}^i_i = u_i$. This is because, since $c_i \leq 0$, the maximum value of $x_i$ in the domain will minimize $c_i x_i$. Moreover, for this choice of $\bar{x}^i_i, \bar{y}^i_i = \frac{r}{u_i}$ is also minimum, and since $d_i \geq 0, \left( u_i, \frac{r}{u_i} \right)$ minimizes $c_i x_i + d_i y_i$.

\textsc{Case 2:} If $c_i > 0$ and $d_i = 0, \bar{x}^i_i = 1, \bar{y}^i_i = r$ as $\bar{x}^i_i \geq 1$.

\textsc{Case 3:} The remaining case is $c_i > 0, d_i > 0$. Since the points $\mathcal{L} \left( i, p_i, \frac{r}{p_i} \right), p_i \in \{ 1, \ldots, u_i \}$ are the extreme points of $conv \left( S^U_i \right)$, minimizing $c_i x_i + d_i y_i$ over $conv \left( S^U_i \right)$ and over $E_i$ are equivalent. To solve this we will use the same analysis as in the proof of Proposition~\ref{2DPol} with slight modification as there is an upper bound $u_i$ on the variable $x_i$. So, in this case we have the following choice of $\bar{x}^i_i$ and consequently $\bar{y}^i_i = \frac{r}{\bar{x}^i_i}$.
\[
\bar{x}^i_i = 
\begin{cases}
 \sqrt{\frac{r d_i}{c_i}}, \ \text{if } \sqrt{\frac{r d_i}{c_i}} \in \{ 1, \ldots, u_i \}, \\
 
 1,  \ \text{if } \sqrt{\frac{r d_i}{c_i}} < 1, \\
 
 \left\lceil \sqrt{\frac{r d_i}{c_i}} \right\rceil \ \text{or } \left\lfloor \sqrt{\frac{r d_i}{c_i}} \right\rfloor, \ \text{whichever minimizes } c_i x_i + d_i \frac{r}{x_i}, \\
 \hspace{1.3in} \text{if } 1 < \sqrt{\frac{r d_i}{c_i}} < u_i \ \text{and } \sqrt{\frac{r d_i}{c_i}} \notin \mathbb{Z}_+,\\
 
 u_i, \ \text{if } \sqrt{\frac{r d_i}{c_i}} > u_i.
\end{cases}
\]

From the above analysis, we can solve the problem (\ref{Pi}) in linear time in the input size, as we just have to check the signs of $n - 1$ entries and have to check the value of $\sqrt{\frac{r d_i}{c_i}}$, whenever it exists, and if not then the signs of $c_i$ and $d_i$. Since finding $\zeta = \min \{ \zeta_i : i \in N \}$ takes $O(n)$ time, we can solve (\ref{P}) in $O(n^2)$ time in the input size.

\end{appendices}


\bibliography{MIBilin} 
\bibliographystyle{plain}


\end{document}